\pdfoutput=1
\documentclass[a4paper,times,3p]{elsarticle}
\usepackage[utf8]{inputenc}

\usepackage[T1]{fontenc}

\usepackage{microtype}

\usepackage{booktabs}

\usepackage{amsmath}
\usepackage{amssymb}
\usepackage{amsthm}
\usepackage{bm}
\usepackage{graphicx}
\usepackage{subfigure}                           
\usepackage{psfrag}
\usepackage{color}

\usepackage{mathtools}

\definecolor{cmykcyan}{cmyk}{1,0,0,0}
\definecolor{cmykred}{cmyk}{0,1,1,0}
\definecolor{cmykblack}{cmyk}{0,0,0,1}

\usepackage[font=footnotesize,labelfont=bf]{caption}
\usepackage{multirow}

\usepackage{pxfonts}

\usepackage{import}

\graphicspath{./graphics/}

\usepackage[textsize=scriptsize]{todonotes}
\usepackage{setspace}

\usepackage{hyperref}

\journal{Computer Methods in Applied Mechanics and Engineering}

\def\statement{\begin{minipage}[t]{.75\textwidth}
       NOTICE: This is the author's version of a work that was accepted for
publication in Computer Methods in Applied Mechanics and Engineering.  Changes
resulting from the publishing process, such as editing, corrections, structural
formatting, and other quality control mechanisms may not be reflected in this
document. Changes may have been made to this work since it was submitted for
publication.
\\

\copyright \, 2017. This manuscript version is made available under the
CC-BY-NC-ND 4.0 license \\
\url{http://creativecommons.org/licenses/by-nc-nd/4.0/}.
       \end{minipage}}

\makeatletter
\def\ps@pprintTitle{%
     \let\@oddhead\@empty
     \let\@evenhead\@empty
     \def\@oddfoot{\footnotesize\itshape
       \statement\hfill\today}%
     \let\@evenfoot\@oddfoot}
\makeatother

\begin{document}
\begin{frontmatter}

\title{Fluid-Structure Interaction with NURBS-Based Coupling}
\author{Norbert Hosters\corref{cor1}}
\ead{hosters@cats.rwth-aachen.de}
\author{Jan Helmig\corref{}}
\ead{helmig@cats.rwth-aachen.de}
\author{Atanas Stavrev\corref{}}
\ead{stavrev@aices.rwth-aachen.de}
\author{Marek Behr\corref{}}
\ead{behr@cats.rwth-aachen.de}
\author{Stefanie Elgeti\corref{}}
\ead{elgeti@cats.rwth-aachen.de}

\cortext[cor1]{Corresponding author}

\address{Chair for Computational Analysis of Technical Systems (CATS), Center for Computational Engineering Science (CCES),\\
	RWTH Aachen University, 52056 Aachen, Germany}

\begin{abstract}
Engineering design via CAD software relies on Non-Uniform Rational B-Splines (NURBS) as a
means for representing and communicating geometry. Therefore, in general, a NURBS
description of a given design can be considered the exact description. The development of
isogeometric methods has made the geometry available to analysis
methods~\cite{Hughes2005}. Isogeometric analysis has been particularly successful in
structural analysis; one reason being the wide-spread use of two-dimensional finite
elements in this field. For fluid dynamics, where three-dimensional analysis is usually
indispensable, isogeometric methods are more complicated, yet of course not impossible, to
apply in a general fashion. This paper describes a method that enables the solution of
fluid-structure-interaction with a matching spline description of the interface. On the structural side, the spline is used in an isogeometric setting. On the fluid side, the same spline is used in the framework of a NURBS-enhanced finite element method (extension of \cite{Sevilla2011a}). The coupling of the structural and the fluid solution is greatly facilitated by the common spline interface. The use of the identical spline representation for both sides permits a direct transfer of the necessary quantities, all the while still allowing an adjusted, individual refinement level for both sides.
\end{abstract}

\end{frontmatter}

\section{Introduction} \label{sec-intro}

\pdfoutput=1

Geometries for engineering applications are generated using Computer-Aided-Design (CAD)
systems: the CAD model is what we assume to be the exact geometry. Nowadays, all major CAD
systems share one common basis for geometry representation: Non-Uniform Rational B-Splines
(NURBS). NURBS provide a common standard to exchange  geometry information. In the classic
discretization methods, the exact NURBS geometry is lost when a finite element mesh ---
which is in general only an  approximation of the geometry --- is generated. Spline-based
methods, such as Isogeometric Analysis \cite{Hughes2005, Cottrell2009} or NURBS-Enhanced
Finite Elements \cite{Sevilla2011a, Stavrev2015} make use of the NURBS format in order to
integrate the exact geometry into the finite element method, thus avoiding the
approximative character of the finite element mesh. Note, the change of computational
domain due to numerical analysis with isogeometric analysis, e.g., for structural
deformation, is only an approximation, even though the NURBS representation of the initial
geometry can be considered as exact.  This paper presents a fluid-structure-interaction (FSI) method with spline-based, smooth geometry representation of the interface between fluid and structure.

To give a general overview, we will categorize different approaches to FSI according to (1) the solution process, (2) the temporal coupling, and (3) the spatial coupling.

(1) With respect to the solution method, one differentiates between monolithic approaches (e.g., \cite{Huebner2004a,Heil2004a,Zilian2010}) or partitioned approaches (e.g., \cite{LeTallec2001,Piperno2001,Foerster2007a,Kuettler2008a,Bazilevs2012a}). The monolithic approach includes the full solution process in one single solver: The main advantage gained is the robustness (in the sense of stability and the capability for large time steps); the price to pay is the significant implementation effort. Partitioned approaches rely on individual solvers for all physical systems, which generate individual contributions to the coupled system. The necessary communication of solution data is restricted to the interfaces between the single field solvers. This facilitates the use of existing, already very evolved solvers, but leads to stability issues. Only in the case of partitioned solvers, there is a need for temporal and spatial coupling. 

(2) Methods for temporal coupling can be categorized into strong and weak coupling schemes. When using a strong coupling method, the idea is to iterate between the individual solvers until a converged solution for the current time step is obtained \cite{Kuettler2008a}. Within a weak coupling, each individual solver is only called once per time step; this usually requires the use of predictor-corrector methods \cite{Piperno2001,Ritter2011}. Weak coupling schemes profit from reduced computational cost. However, they may suffer from a lack of stability in cases where the density ratio between fluid and solid is close to one \cite{Joosten2010a}. This limit may very well be reached for incompressible flows in combination with flexible structures.

(3) With regard to the spatial coupling, the use of individual solvers leads to \--- in general \--- non-conforming discretizations for fluid and structure. Consequently,  the transfer of interfacial quantities (e.g., forces, displacements) from one computational mesh to the other requires projection methods \cite{Felippa2010}. Again, there are two fundamental approaches available: mesh-based and mesh-free methods. Mesh-based methods, such as finite-interpolation-elements (FIE), rely on the direct approaches of the individual solvers \cite{Beckert2000,Kollmannsberger2007}. An example for a mesh-free method is the spline-based moving-least-squares-method (MLS) \cite{Quaranta2005}, which interpolates the transfer information between the different interface discretizations. Both methods can also be adapted to individual solvers with spline-based solution methods \cite{Bazilevs2008a,Kloeppel2011,Farhat1998}. 

In this work, we use a partitioned approach with strong coupling and rather specialized
solvers. The key idea of our FSI
approach is to generate a smooth and conforming spline description of the interface
between fluid and structure, that might be also benefical for monolithic schemes. On the fluid side, the spline is supported through the NURBS-Enhanced Finite Element Method. On the structural side, this is done by means of IGA. The use of the identical spline representation for both sides permits a direct transfer of the necessary quantities, all the while still allowing an adjusted, individual refinement level for both sides. Although there are several approaches that use IGA on the structural side, to our knowledge this is the first time where a conforming interface for both sides was established using splines. 

\begin{figure}[h]
\footnotesize
\subfigure[FSI interface representation with piecewise linear discretization.  ]{
	\centering
	\includegraphics[width=6.5cm]{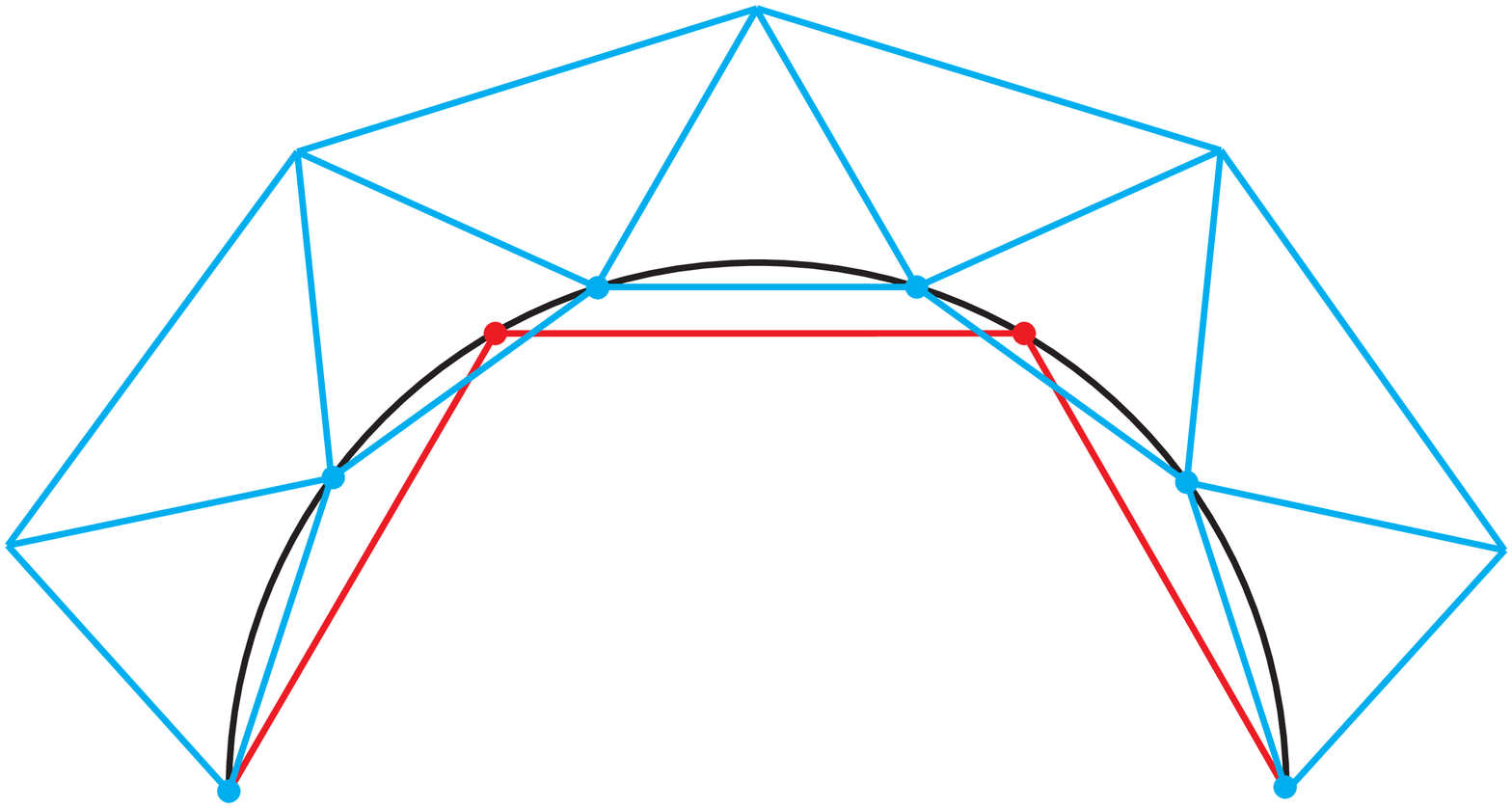}
	\label{CurrentSituation}
}
\hspace{0.5cm}
\begin{minipage}[h]{1cm}
	\centering
	\includegraphics[width=1cm]{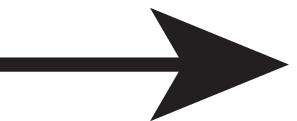}
\end{minipage}
\hspace{0.5cm}
\subfigure[FSI interface representation using identical splines.]{
	\centering
	\includegraphics[width=6.5cm]{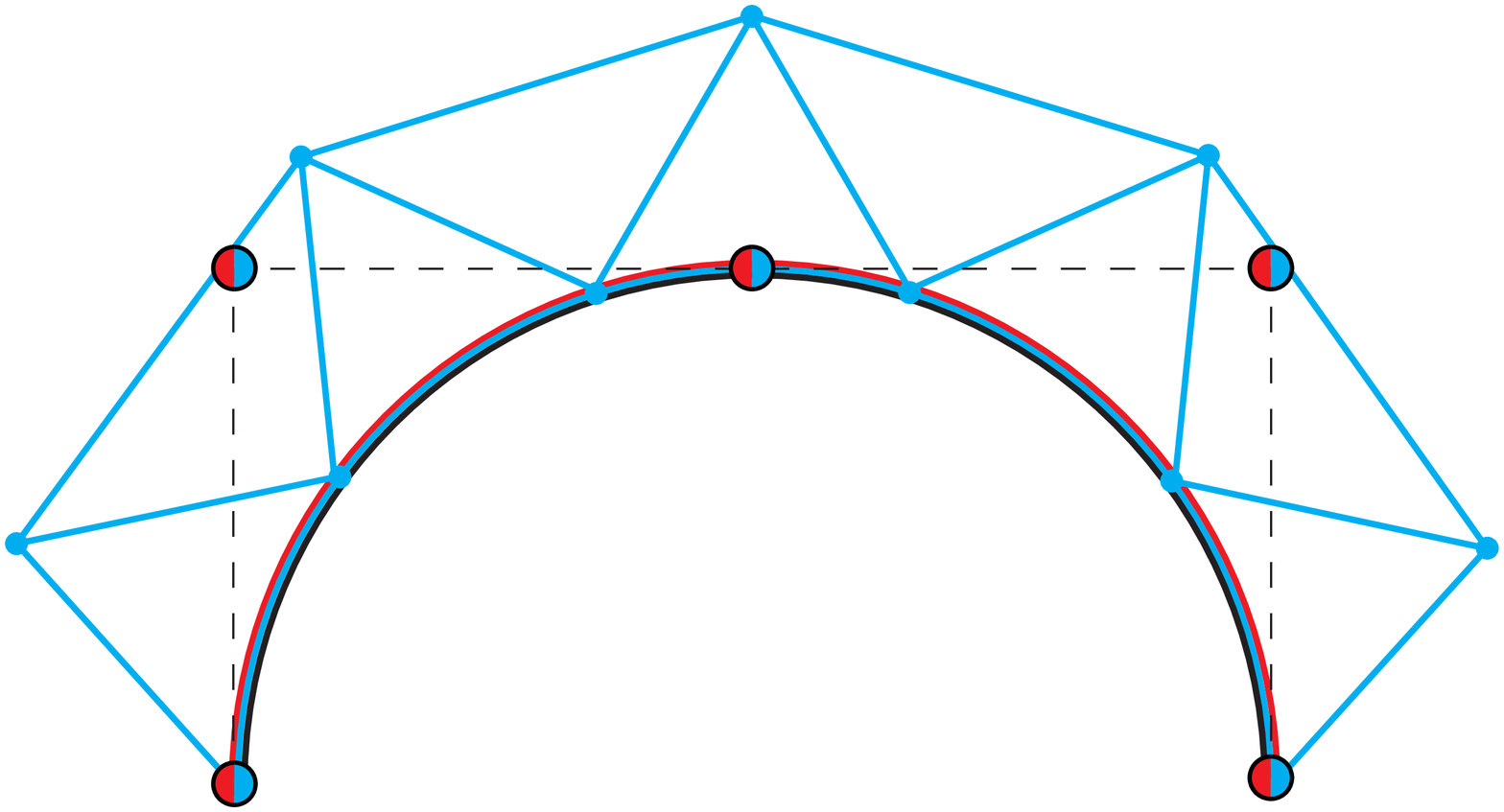}
	\label{ACM2IGA}
}
\newline
\begin{center}
\begin{minipage}[]{0.75\linewidth}
{\color{cmykcyan}---} CFD mesh \hfill\includegraphics[width=0.2cm]{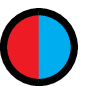} Spline control points\\
{\color{cmykblack}---} Actual interface \\
{\color{cmykred}---} Structural mesh
\end{minipage}
\end{center}
\caption[]{\footnotesize Discretization of a given, curved interface (black). In (a), we see \--- for both the fluid and the structure \--- a piecewise linear interface discretization. Note that the two discretizations of the same interface are incompatible with each other, making projection inevitable for the transferred quantities (forces, displacements, etc.). In particular for curved interfaces, this necessarily results in numerical errors, which contribute to the destabilization of the partitioned approach. In (b), the interface is now represented using the identical spline curve for both sides. The conforming and smooth interface representation reduces the numerical error to the discretization error even for partitioned FSI approaches. }
\end{figure}

The paper has five main sections. Section 2 is devoted to the underlying governing equations. Sections 3--5 outline the numerical methods for fluid, structure, and coupling. Numerical examples are discussed in Section 6.

\section{Modeling: Governing Equations} \label{sec-governing}

The test cases in this work feature a nonlinear elastic structure enclosed in an incompressible Newtonian fluid. This section introduces the relevant governing equations for the structure and the fluid.

\subsection{Structural Deformation} 

The response of a structure to an external load is governed by classical elastodynamics. The primal variable is the displacement field, here denoted as ${\bf{d}}^s \left( {\bf{x}}, t \right)$. It describes the change between the current structural configuration ${\bf{x}}$ and the initial configuration ${\bf{x}}_0$. These alterations can be either due to deformations of the structure or due to rigid body motions. 

Consider a deformable computational domain, which is, at each instant in time, denoted by \(\Omega_t^s \). \(\Omega_t^s \) is a subset of \(\mathbb{R}^{nsd^s}\), with \(nsd^s\) as the number of space dimensions for the structure. Its boundary is denoted as  $\Gamma _t ^s$. At each point in time \(t \in [0,T] \), the displacement ${\bf{d}} \left( {\bf{x}}, t \right)$ of the structure is governed by an equation based on Newton's second law:

\begin{equation}\label{Elastodynamics}
 \rho ^s \frac{d^2 {\bf{d}}^s}{d t ^2} = \boldsymbol{\nabla} \cdot \boldsymbol{\sigma} ^s + {\bf{b}} ^s \quad \mbox{on} \  (\Omega_t^s) \ \forall t \in [0,T] \,,
\end{equation}

\noindent where $\rho ^s$ indicates the density of the structure, ${\bf{b}}^s $ the prescribed body force per unit volume and $\boldsymbol{\sigma ^s}$ the Cauchy stress tensor. 

As constitutive relation, we employ a geometrically nonlinear approach for a hyperelastic material model; a model which enables us to account for large deformations. More specifically, we employ the St. Venant Kirchhoff model, which describes the constitutive equation for the stress tensor based on a stress-strain relation using the Green-Lagrange strain tensor ${\bf{E}}^s$ and the 2nd Piola-Kirchhoff stress tensor $\bf{S}^s(\bf{E}) $ \cite{Fu2001}. They are defined as:

\begin{equation}
{\bf{S}}^s(\bf{E})=  \lambda^s tr\left( {\bf{E}}^s \right) + 2 \mu^s {\bf{E}^s}
\end{equation}

\noindent and 

\begin{equation}
{\bf{E}} = \frac{1}{2} \left( {\bf{F}}^T {\bf{F}} - {\bf{I}} \right),
\end{equation}

\noindent where the deformation gradient ${\bf{F}}$ is defined as

$$ \bf{F} = \frac{\partial{\bf{x}}}{\partial{\bf{x}_0}}. $$

Here, $\lambda^s$ and $\mu^s$ are the Lam\'e parameters. They can be expressed in terms of the Young's modulus $E^s$ and the Poisson ratio $\nu^s$ as

\begin{equation}
\lambda^s = \frac{\nu^s E^s}{(1+\nu^s)(1-2\nu^s)}, \quad
\mu^s = \frac{E^s}{2(1+\nu^s)}.
\end{equation}

Connecting all of this information, we can express the equation of motion in the reference configuration as

\begin{equation}\label{Elastodynamicsnew}
 \rho ^s \frac{d^2 {\bf{d}}^s}{d t ^2} = \boldsymbol{\nabla}_0 \cdot \left( {\bf{S}} ^s {\bf{F}}^T \right) + {\bf{b}} ^s \quad \mbox{on} \  (\Omega_0^s) \ \forall t \in [0,T] \,.
\end{equation}

In order to obtain a well-posed system, boundary  conditions have to be imposed on $\Gamma^s$. Here, we distinguish between Dirichlet $({\bf{g}}^s )$ and Neumann  $({\bf{h}}^s )$ boundary conditions given by:

\begin{equation}\label{bcstructure}
  {\bf{d}}^s = {\bf{g}}^s \; \mathrm{on} \; \left( \Gamma ^s_t \right) _g ,  \quad {\bf{n}}^s \cdot \boldsymbol{\sigma}^s = {\bf{h}}^s \; \mathrm{on} \; \left( \Gamma ^s_t \right) _h.
\end{equation}

$ \left(\Gamma_t^s\right)_g $ and $ \left(\Gamma_t^s\right)_h $ denote the Dirichlet and Neumann part of the boundary, forming a complementary subset of $\Gamma^s_t$, i.e., $ \left(\Gamma_t^s\right)_g \cup \left(\Gamma_s^f\right)_h = \Gamma^s_t$ and $  \left(\Gamma_t^s\right)_g \cap \left(\Gamma_t^s\right)_h = \emptyset$. Here, ${\bf n}^s$ refers to the outer normal vector on $ \Gamma_t^s $. 

Furthermore, the initial displacement is prescribed as initial condition:

\begin{equation}\label{initcondelastodynamics}
{\bf{d}}^s \left( {\bf{x}}, t = 0 \right) =  {\bf{d}}^0.
\end{equation}

\subsection{Fluid Flow}

Consider a deformable fluid domain, which is, at each instant in time, denoted by \(\Omega_t^f \). Domain \(\Omega_t^f \) is a subset of \(\mathbb{R}^{nsd^f}\), with \(nsd^f\) as the number of space dimensions of the fluid domain. Then at each point in time \(t \in [0,T] \), the velocity, \( {\bf u}^f({\bf x},t) \), and the pressure, \( p^f({\bf x},t) \), of the fluid are governed by the unsteady, incompressible Navier-Stokes equations:

\begin{align}
\label{NavierStokes1}
\rho^f \left( \frac{\partial {\bf u}^f}{\partial t} + {\bf u}^f \cdot \nabla {\bf u}^f - {\bf f}^f \right) -  \nabla \cdot \bm{\sigma}^f &= {\bf 0} \quad \mbox{on} \ (\Omega_t^f) \ \forall \  t \in [0,T] \,, \\
\label{NavierStokes2}
{\bf \nabla \cdot u }^f &= 0 \quad \mbox{on} \  (\Omega_t^f) \ \forall t \in [0,T] \,,
\end{align}
with \( \rho^f\) as the fluid density. In the Newtonian case, the stress tensor \( \bm{\sigma}^f\) is defined as 

\begin{align}
\bm{\sigma}^f ( {\bf u}^f, p^f) &= -p^f {\bf I} + 2 \mu \bm{\varepsilon}^f( {\bf u}^f) \quad \mbox{on} \quad (\Omega_t^f) \,, 
\label{eqn:Newtonian}
\end{align}
with
\begin{align} \label{eqn:strain}
\bm{\varepsilon}^f({\bf u}^f) &= \frac{1}{2} \left(\nabla {\bf u}^f + (\nabla {\bf u}^f)^T\right) \,,
\end{align}
where \( \mu^f  \) denotes the dynamic viscosity.  $ {\bf f}^f $ includes all external body forces per unit mass of fluid. 

In order to obtain a well-posed system, boundary  conditions have to be imposed on the external boundary of $\Omega^f_t$, denoted as $\Gamma^f_t$. Here, we distinguish between Dirichlet and Neumann boundary conditions given by:

\begin{align}
{\bf u}^f &= {\bf g}^f \quad \mbox{on} \ \left(\Gamma_t^f\right)_g, \\
{\bf n}^f \cdot \bm{\sigma}^f &= {\bf h}^f \quad \mbox{on} \ \left(\Gamma_t^f\right)_h \label{eqn:Neumanncondition} \,,
\end{align} 

\noindent where ${\bf g}^f$ and ${\bf h}^f$ are prescribed velocity and stress values. $ \left(\Gamma_t^f\right)_g $ and $ \left(\Gamma_t^f\right)_h $ denote the Dirichlet and Neumann part of the boundary, forming a complementary subset of $\Gamma^f_t$, i.e., $ \left(\Gamma_t^f\right)_g \cup \left(\Gamma_t^f\right)_h = \Gamma^f_t$ and $  \left(\Gamma_t^f\right)_g \cap \left(\Gamma_t^f\right)_h = \emptyset$. Here, ${\bf n}^f$ refers to the outer normal vector on $ \Gamma_t^f $. 

In the transient case, a divergence-free velocity field for the whole computational domain is needed as an initial condition:

\begin{align}
{\bf u}^f ( {\bf x}, 0 ) = {\bf u}^0( {\bf x}) \quad \mbox{in}\ \Omega_t^f \ \mbox{at} \  t=0 \,.
\end{align}

\subsection{Coupling Conditions at the Fluid-Structure Interface}

The governing equations of fluid (Equations~\eqref{NavierStokes1}--~\eqref{NavierStokes2})
and structure (Equation~\eqref{Elastodynamicsnew}) need to be connected in order to
represent the interaction between the two components. This interaction takes place only
through the common interface $\Gamma_{FS} = \Gamma^{f}_t \cup \Gamma^{s}_t$, not the full volume; thus distinguishing FSI problems from other multiphysics problems \cite{Wall99}.

For a consistent coupling, the following physical requirements are essential: (1)
geometric compatibility between the fields, (2) kinematic and dynamic conditions at the
shared interface $\Gamma_{FS}$, and (3) conservation of mass, momentum and energy. This
leads to the following coupling conditions on $\Gamma_{FS}$:

\noindent \emph{Kinematic continuity:}

\begin{equation}\label{eq:kinematic}
 \begin{split}
 {\bf d}^f \left( {\bf x}, t \right) &= {\bf d} ^s  \left(  {\bf x}, t \right) \quad\quad   \mathrm{on} \quad \Gamma_{FS}, \\
 {\bf u}^f \left( {\bf x}, t \right) &= {\bf u} ^s \left(  {\bf x}, t \right)  \quad \quad \mathrm{on} \quad \Gamma_{FS}.
 \end{split}
\end{equation}

\noindent These coupling conditions ensure the continuity of displacements and velocities across the interface. 

\noindent \emph{Dynamic continuity:}

\begin{equation}\label{eq:dynamic}
\boldsymbol{\sigma} ^ f \left(  {\bf{x}}, t \right) \cdot {\bf{n}}^f = - \boldsymbol{\sigma} ^s \left(  {\bf{x}}, t \right) \cdot {\bf{n}}^s \quad \quad \mathrm{on} \quad \Gamma_{FS}.
\end{equation}

\noindent In agreement with Newton's third law \--- $Actio \; and \; Reactio$  \--- this coupling condition enforces continuity of fluid ($\boldsymbol{\sigma}^f$) and structural stresses ($\boldsymbol{\sigma}^s$) at the interface $\Gamma_{FS}$  \cite{Braun2007}. 

\section{Numerical Methods: Structural Solution} \label{sec-structuralsolver}

The governing equation for the structural deformation, Equation~\eqref{Elastodynamicsnew},
is discretized using an isogeometric finite element method. The time discretization is
performed using a generalized $\alpha$ scheme \cite{Kuhl1999, chung_time_1993}. The resulting weak form is linearized using a Newton-Raphson framework.

\subsection{Variational Form}

The employed isogeometric finite element method is a standard Galerkin formulation; only
the interpolation functions differ from standard Lagrange interpolation functions. We
introduce a finite-dimensional space $\mathcal{I}_0 \subset C^0(\overline{\Omega^s})$ based on
NURBS interpolation as spatial discretization. To obtain the weak form, Equation~\eqref{Elastodynamicsnew} is multiplied by the displacement test function ${\bf w}^s$, integrated, and the stress term $\bm{\sigma}^s$ is integrated by parts. The following finite element interpolation and weighting function spaces for the displacement ${\bf d}^s$ can be defined: 

\begin{align} \label{eq:function_spaces_structure}
\tilde{\mathcal{S}} _d ^{s,h} &= \{ \; {\bf d}^{s,h} | \; {\bf d}^{s,h} \in
\mathcal{I}_0,\, {\bf d}^{s,h} \doteq 
  {\bf g}^{s,h}  \quad \mathrm{on} \quad \Gamma ^s _g,\,\forall t \in [0,T]\;  \}, \\
 \tilde{\mathcal{V}} _d ^h &= \{ \; {\bf w}^{s,h} | \; {\bf w}^{s,h} \in\mathcal{I}_0,\,
{\bf w}^{s,h} \doteq {\bf 0}  \quad  \mathrm{on} \quad  \Gamma ^s _g,\,\forall t \in [0,T] \; \}.
 \end{align}

Based on these spaces, the discrete formulation of the structure is defined as: Find $ {\bf d}^{s,h} \in \tilde{\mathcal{S}} _d ^h$ such that:

\begin{equation}\label{weakstructure}
\delta W \left( {\bf w}^{s,h}, {\bf d}^{s,h} \right) = \int_{\Omega ^s} {\bf w}^{s,h}
\cdot \frac{d^2 {\bf d}^{s,h}}{d t ^2} +  \bm{\nabla} {\bf w}^{s,h} \colon \bm{\sigma}^s \; d \Omega - \int_{\Gamma_g ^s } {\bf w}^{s,h} \cdot {\bf h}^{s,h} \; d \Gamma = 0
\end{equation}

\noindent holds for all ${\bf w}^{s,h} \in \tilde{\mathcal{V}} _d ^h$. \\

\subsection{Isogeometric Analysis (IGA)}

Both the discrete function spaces (Equation~\eqref{eq:function_spaces_structure}) and the weak form (Equation~\eqref{weakstructure}) assume the use of isogeometric interpolation functions. In our implementation, we have resorted to a standard IGA method, as it was introduced in \cite{Hughes2005, Cottrell2009}. In view of the ample publications in the field of structural analysis with isogeometric finite elements, e.g., \cite{Dornisch2013}, we will restrict ourselves to describing the general idea of IGA. 
IGA invokes an isoparametric finite element concept. This means that the unknown solution and the geometry are interpolated utilizing the same type of interpolation functions. In the case of IGA, these interpolation functions are inspired by the standard CAD description employed in design engineering. CAD systems are usually based on either NURBS \cite{Piegl97} or, in very rare cases, also on T-splines \cite{Bazilevs2010}. Both cases belong to the category of parametric geometry description, with basis functions defined in terms of a local parameter and a mapping to the physical space; a concept very similar to Lagrange interpolation functions. 
The key advantages of employing the NURBS or T-spline basis are (1) exact geometry
description in the initial state and (2) user-controlled smoothness of the basis.

A NURBS curve is represented through a combination of control points that guide the curve --- ${\bf P}_i$ indicates the control point coordinates --- and the NURBS basis functions $R_{i,p}$ of degree $p$. The curve definition is as follows:

\begin{align}
{\bf C}(\Theta) = \sum_{i=1}^n R_{i,p}(\Theta) {\bf P}_i,
\end{align}

\noindent with $\Theta$ as the local parameter of the NURBS curve. Due to the isoparametric concept, the unknown solution ${\bf d}^{s,h}(\Theta)$ is represented in the following form:

\begin{align}
{\bf d}^{s,h}(\Theta) = \sum_{i=1}^n R_{i,p}(\Theta) {\bf d}^s_i.
\end{align}

\begin{figure}
\centering
\psfrag{a}{$u_1$}
\psfrag{b}{$u_2$}
\psfrag{c}{$u_3$}
\psfrag{d}{$u_4$}
\psfrag{e}{$u^h$}
\psfrag{f}{$x$}
\psfrag{g}{$u^h(x)$}
\psfrag{h}{${\bf P}_1$}
\psfrag{i}{${\bf P}_2$}
\psfrag{j}{${\bf P}_3$}
\psfrag{k}{${\bf P}_4$}
\psfrag{l}{$\Theta$}
\psfrag{x}{$x$}
\psfrag{p}{$0.0$}
\psfrag{q}{$R_i(\Theta)$}
\psfrag{r}{$1.0$}
\psfrag{m}{$0.0$}
\psfrag{n}{$\Theta$}
\psfrag{o}{$1.0$}
\includegraphics[scale=0.45]{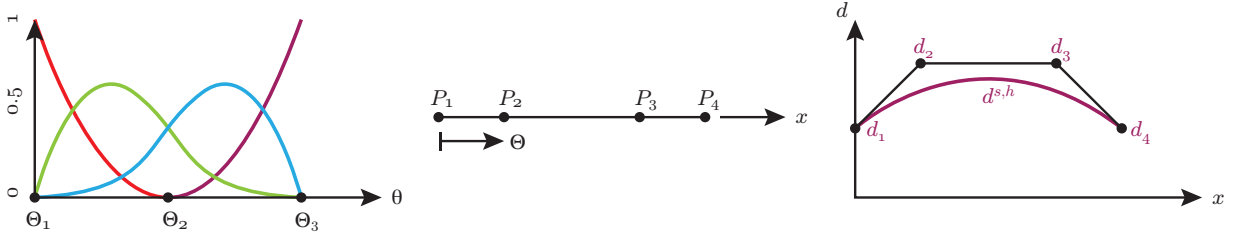}
	\caption{Geometry and function representation in IGA: In the left picture, sample basis functions for a quadratic NURBS are illustrated. On the one hand this basis can be employed to interpolate control points in order to represent a geometry (e.g., the rod in the middle picture). On the other hand, it can be used to interpolate the unknown function $d^{s,h}$. In this case, the discrete solution is represented through so-called control variables  $d_1$ through $d_4$.  Note that these control variables do not necessarily coincide with the solution curve; a property connected to the structure of the basis functions $R_i$. The local parameter along the curve is $\Theta$, and it is embedded into $\mathbb{R}^1$ with coordinate $x$. }
\label{fig:iga-function-representation}
\end{figure}

\section{Numerical Methods: Flow Solution} \label{sec-flowsolver}

The Navier-Stokes equations \eqref{NavierStokes1}--\eqref{NavierStokes2} are discretized using P1P1 finite elements, i.e., linear interpolation for both the velocity and pressure degrees of freedom \cite{Tezduyar92a,Behr92d}. Along the boundary, the linear finite elements are supplemented by geometry information using the NURBS-Enhanced Finite Element Method described in Section \ref{sec-NEFEM} \cite{Stavrev2015}. P1P1 finite elements are known to violate the Ladyzhenskaya-Babuska-Brezzi (LBB) compatibility condition. Consequently, without appropriate stabilization, the pressure field is likely to present spurious and oscillatory results. The stabilization technique used here is Galerkin/Least-Squares (GLS) stabilization. In the GLS method, the stabilization term consists of an element-by-element weighted least-squares form of the original differential equation \cite{Donea2003}.

\begin{figure}[htbp]
\center
\includegraphics[width=8.0cm]{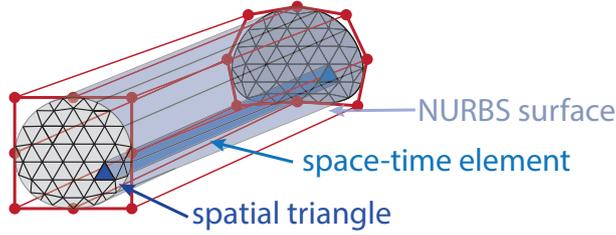}
\caption{Illustration of the space-time method on deforming domains: The spatial domain is a 2D-circle, which deforms over one timestep. The spatial domain is meshed with triangles. As it is extruded into the time direction, it becomes a space-time prism. Through the use of NEFEM, the domain boundary is represented using a non-uniform rational B-spline (indicated in the picture through the red control points). Elements on the boundary have one curved edge defined by the spline. The finite element nodes coincide with specific points on the spline. Therefore, each finite element node is also assigned a spline coordinate $\Theta_{FE-node number}$.}
\label{fig:spacetimeNEFEM}
\end{figure}

The required deformability of the computational domain over time is included through the Deformable-Spatial-Domain/Stabilized Space-Time (DSD/SST) method~\cite{Tezduyar92a}. In contrast to the classic choice within the finite element community of finite element discretization in the spatial domain and finite differences in the time direction, this method relies on finite elements in both space and time (cf. Figure~\ref{fig:spacetimeNEFEM}). This approach allows to formulate the variational form directly over the deforming domain. The necessity of modifying the fluid-flow equations based on the mesh deformation velocity is avoided. As long as the domain deformation remains within a certain range, the mesh can be equipped with the necessary flexibility to account for the mesh motion without the need for remeshing.  For this purpose, the Elastic Mesh Update Method~\cite{Johnson94a} is employed. 

\subsection{Variational Form} \label{sec-variationalform}

In order to construct the finite element function spaces for the space-time method, the time interval $[0,T]$ is divided into subintervals \(I^f_n = (t_n,t_{n+1})\), with \(t_n\) and \(t_{n+1}\) representing an ordered series of time levels \( 0=t_0 < t_1< \cdots < t_N =T\).  Now, if \(\Omega_n^f=\Omega_{t_n}^f \), the space-time slab \(Q_n^f \) is defined as the domain enclosed by the surfaces \( \Omega_n^f, \Omega_{n+1}^f\) as well as the surface desribed by $\partial \Omega_t^f$ as t traverses \(I_n^f\), which shall be named \(P_n^f\). The following finite element interpolation and weighting function spaces for velocity ${\bf u}$ and pressure $p$ can be defined based on linear, $C^0$-continuous interpolation in space and also linear, but discontinuous interpolation in time: 

\begin{align}
(\mathcal{S}^h_{\bf u})_n &= \{ {\bf u}^{f,h} | {\bf u}^h \in [H^{1h}(Q_n)]^{n_{sd}} , {\bf u }^{f,h} \doteq {\bf g}^{f} \quad \text{on} \quad (P_n)_{\bf g} \}, \\
(\mathcal{V}^h_{\bf u})_n &= \{ {\bf w}^{f,h} | {\bf w}^{f,h} \in [H^{1h}(Q_n)]^{n_{sd}} , {\bf w }^{f,h} \doteq {\bf 0} \quad \text{on} \quad (P_n)_{\bf g} \}, \\
(\mathcal{S}^h_p)_n &= (\mathcal{V}^h_p)_n = \{ p^h|p^h \in H^{1h}(Q_n)\}.
\end{align}

The stabilized space-time formulation of the incompressible Navier-Stokes equations (\ref{NavierStokes1})--(\ref{NavierStokes2}) for deforming domains can then be expressed as follows:  Given \( ( {\bf u}^h)^-_n\) find \( {\bf u}^h \in (\mathcal{S}^h_{\bf u})_n \) and \(p^h \in  (\mathcal{S}^h_p)_n \) such that \( \forall {\bf w}^h \in (\mathcal{V}^h_{\bf w})_n \),  \( \forall  q^h \in (\mathcal{V}^h_p)_n \):
\begin{align}
\label{Variation}
\begin{split}
\int_{Q_n} {\bf w}^{f,h} \cdot \rho^f \left( \frac{\partial {\bf u}^{f,h}}{\partial t} + {{\bf u}^{f,h} \cdot {\bf \nabla u}^{f,h}} - {\bf f}^{f,h} \right)\;dQ + \int_{Q_n} \mbox{\boldmath$\varepsilon$}^{f}({\bf w}^{f,h}):\mbox{\boldmath$\sigma$}^{f,h} (p^{f,h},{\bf u}^{f,h})\;dQ \\
+ \int_{Q_n} q^{f,h}\nabla \cdot {\bf u}^{f,h}\;dQ + \int_{\Omega_n} ({\bf w}^{f,h})^+_n \cdot \rho^f \left( ({\bf u}^{f,h})^+_n - ({\bf u}^{f,h})^-_n \right)\;d\Omega\\
+ \sum_{e=1}^{({n_{el})}_n} \int_{Q^e_n} \tau_{\mbox{\tiny{MOM}}} \frac{1}{\rho^f} \left[ \rho^f \left( \frac{\partial {\bf w}^{f,h}}{\partial t}+{{\bf u}^{f,h} \cdot {\bf \nabla w}^{f,h}} \right) - \nabla \cdot \mbox{\boldmath$\sigma$}^{f,h}(q^{f,h},{\bf w}^{f,h}) \right] \\
\cdot \left[ \rho^f \left(\frac{\partial {\bf u}^{f,h}}{\partial t} + {{\bf u}^{f,h} \cdot {\bf \nabla u}^{f,h}} - {\bf f}^{f,h} \right) - \nabla \cdot {\boldsymbol{\sigma}}^{f,h}(p^{f,h}, {\bf u}^{f,h}) \right]\;dQ
+ \sum_{e=1}^{({n_{el})}_n} \int_{Q^e_n} \tau_{\mbox{\tiny{CONT}}} \nabla \cdot {\bf w}^{f,h} \rho^{f} \nabla \cdot {\bf u}^{f,h}\;dQ \\ 
= \int_{P_n} {\bf w}^{f,h} \cdot {\bf h}^f \;dP .
\end{split}
\end{align}

In the above equation, the following notation is used:
\begin{align}
({\bf u}^h)^{\pm}_n &= \lim_{\epsilon \to 0} {\bf u}(t_n \pm \epsilon),\\
\int_{Q_n} \dots dQ &= \int_{I_n} \int_{\Omega^h_t} \dots d \Omega dt,\\
\int_{P_n} \dots dP &= \int_{I_n} \int_{\partial \Gamma^h_t} \dots d \Gamma dt.
\end{align}

The problem is solved sequentially for each space-time slab, starting with:
\begin{align}
({\bf u}^h)^+_0 = {\bf u}_0.
\end{align}

Details on the method and its parameters \( \tau_{\mbox{\tiny{MOM}}} \) and \( \tau_{\mbox{\tiny{CONT}}} \) can be found in \cite{Pauli2017}.  

Within the stabilization terms, $\nabla \cdot {\boldsymbol{\sigma}}^{f}$ needs to be
computed. This involves the computation of second order derivatives of the velocity field,
which are zero in case of linear shape functions. Thus, a least-squares recovery technique is
applied to improve consistency of the method~\cite{jansen_better_1999}. By that, the
fluid stresses are computed,
which are needed to evaluate tractions on the surface later on.

\subsection{NURBS-Enhanced Finite Element Method} \label{sec-NEFEM}

Comparable to the IGA approach pursued for the structure (cf. Section~\ref{sec-structuralsolver}), also the fluid side benefits from an integration of the CAD geometry. As an additional aspect, the use of the CAD geometry also in the flow solver ensures an excellent compatibility in the fluid-structure-coupling. In principle, it would be possible to also solve the fluid equation \eqref{Variation} using IGA. However, despite recent progress in a variety of directions \cite{Zhang2012,Schillinger2012}, the generation of closed volume splines describing complex geometries remains a challenge. 
As an alternative on middle ground between isogeometric analysis and standard finite elements, Sevilla, Fernandez-Mendes, and Huerta have proposed the NURBS-Enhanced Finite Element Method (NEFEM) \cite{Sevilla2008a,Sevilla2011e,Sevilla2011a}. Already in \cite{Stavrev2015}, the NEFEM has been extended to space-time finite elements. The FSI-coupling requires further modifications detailed in Section~\ref{NEFEM-new}

\subsubsection{The Original Formulation of NEFEM}

\begin{figure}[htbp]
\center
\subfigure[Position of quadrature points]{
\includegraphics[height=4.0cm]{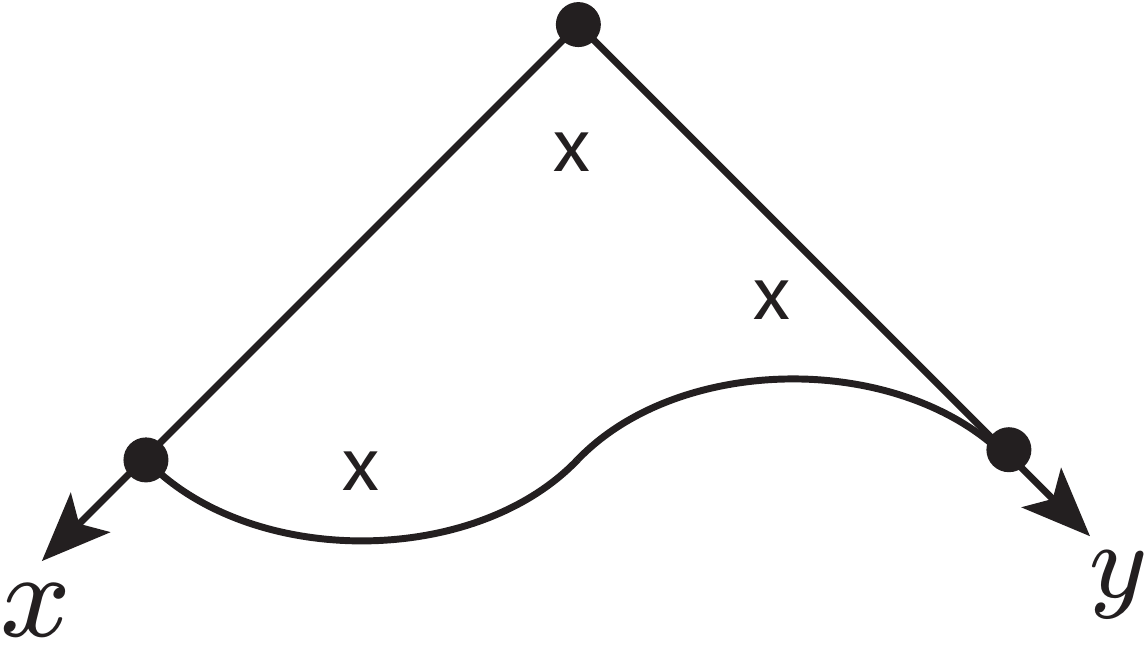}
}
\hfil
\subfigure[Example of a linear shape function]{
\label{resistance}
\includegraphics[trim=4cm 8cm 5cm 10cm, clip, height=5.0cm]{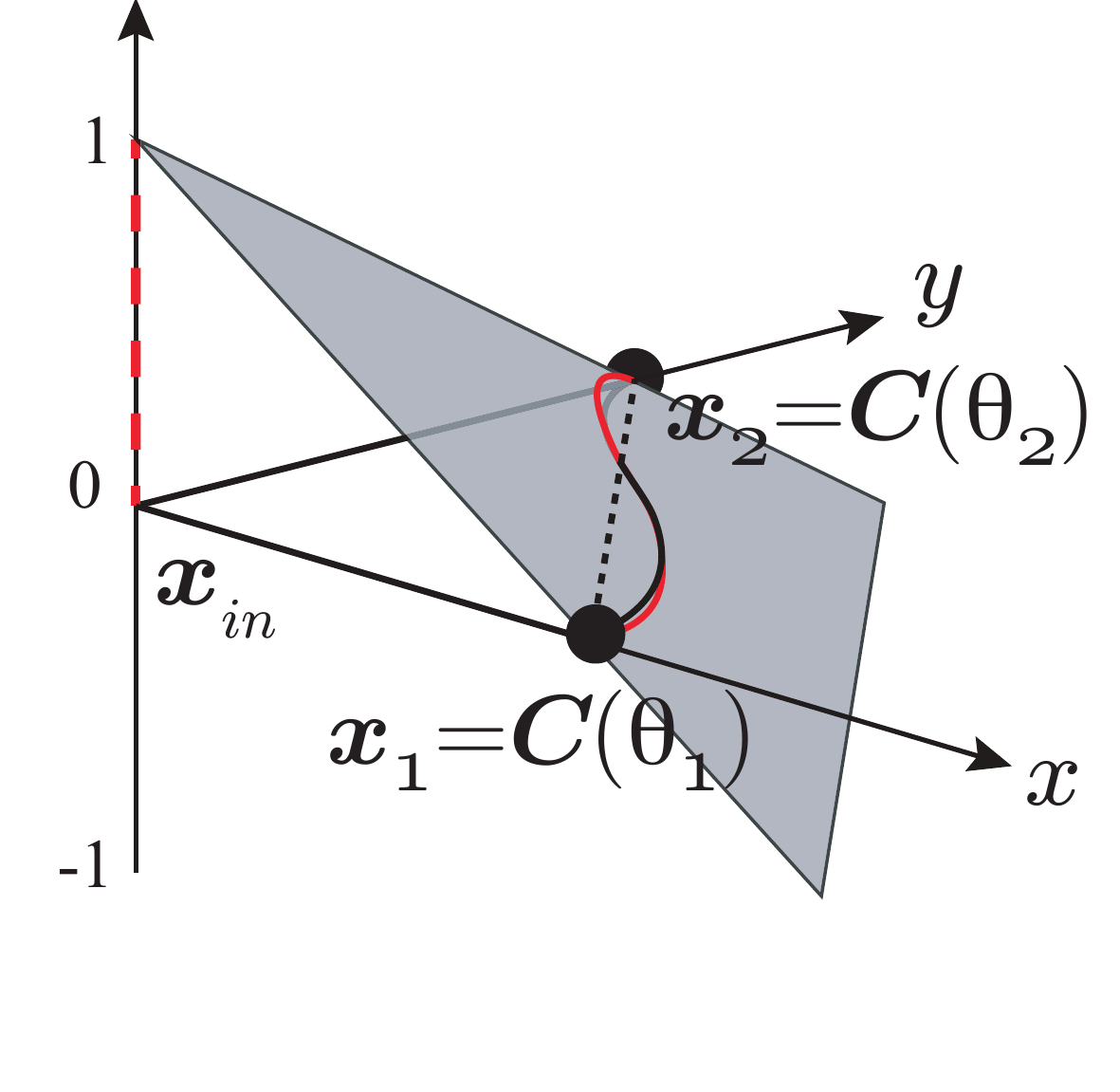}
}
\caption{Illustration of the principal NEFEM concepts as introduced in \cite{Sevilla2008a}. Triangular elements along the boundary of the domain are equipped with one, possibly curved edge, which is represented by a portion of a NURBS curve. This gives access to the geometry. (a) The NEFEM quadrature points x are adapted to the curved triangle shape. (b) Example of a linear shape function in the NEFEM context. Negative values and values larger than $1.0$ may occur.}
\label{fig:NEFEM}
\end{figure}

In NEFEM, instead of the full geometry, only its boundary is represented using NURBS
at the cost of maximally a 2D spline (the boundary of a three-dimensional object). In the interior of the geometry, a standard finite element mesh is utilized, which preserves all advantages of existing meshing algorithms (cf. Figure~\ref{fig:spacetimeNEFEM}). Note that this approach leads to two different kinds of elements: (1) standard finite elements in the interior and (2) elements with a NURBS edge alongside the boundary. Usually, elements of category (1) will be in the vast majority, keeping the computation very efficient. Through the elements of category (2), the geometry is made available during the process of evaluating the integrals of the appropriate weak form (e.g., of the variational equation~\eqref{Variation}): The integration domain $\left(\Omega^f \right)^h$ is no longer only an approximation of the real domain $\Omega^f$. Translated to the numerical implementation, this factor is incorporated through the position of the integration points needed in the finite element method during numerical quadrature: These are determined from the curved NURBS geometry and not from the approximated geometry (cf. Figure~\ref{fig:NEFEM} for a schematic illustration). If we consider boundary integrals --- e.g., for the evaluation of the forces the fluid exerts on the structure --- the integration points are distributed on the curved geometry. 

The unknown function continues to be represented using Lagrange polynomials. These can be of any order, in particular independent of the order of the NURBS basis. In the NEFEM version proposed in \cite{Sevilla2011a,Sevilla2008a,Sevilla2011e} the choice was made to compute the shape functions in the global coordinates. This variant --- termed by the authors as Cartesian FEM --- has the advantage that one obtains undistorted polynomial functions in the global elements, an attribute of particular importance when used with polynomials of higher order. The price to pay is an oddity that arises from this representation: negative shape function values as well as values larger than one may occur along the boundary (cf. Figure~\ref{fig:NEFEM}).

In terms of implementation, there are two important modifications with respect to a standard finite element method: (1) the definition and evaluation of the shape functions and (2) the placement of the quadrature points. Regarding (1), Sevilla performs all evaluations in the global coordinate system. With respect to (2), the original method uses a bi-unit square reference element; even for triangular elements. This enables the use of tensorized 1D Gauss-Legendre quadrature rules. Furthermore, a clear distinction between the NURBS direction and the standard direction can be made; leading to straight interior edges.

\subsubsection{The Modified Version of NEFEM Used in the FSI Context} \label{NEFEM-new}

The NEFEM implementation utilized in this work is based on the same idea as it was proposed by \cite{Sevilla2008a}. The differences lie in the fact that (1) we employ the space-time version, as derived in \cite{Stavrev2015}, and (2) we refrain from using Cartesian FEM, but compute both the shape functions and the position of the integration points on a --- now triangular --- reference element. What remains unchanged is the overall concept of adjusting of the position of the integration points to the curved NURBS shape as well as the representation of the unknown solution with Lagrange polynomials --- even if restricted to linear polynomials in our case. 

\begin{figure}[h]
\centering
 \includegraphics[width=8cm]{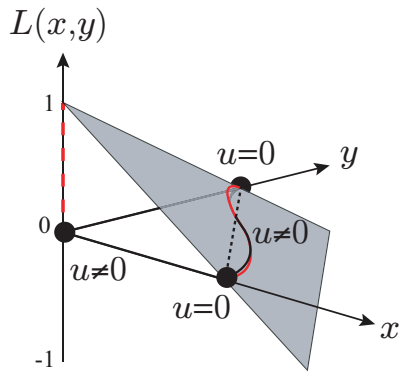}
 \caption{Dirichlet boundary conditions for Cartesian space-time NEFEM. Even though the boundary conditions are set at the nodes of the triangle, the value of the shape function belonging to the interior node is non-zero at the quadrature points, thus leading to an influence of the interior nodal value on boundary integrals.}
 \label{cartesianbc}
\end{figure}

 \begin{figure}[h!]
\centering
  \subfigure[Shape function for a boundary node]{
    \includegraphics[width=6cm]{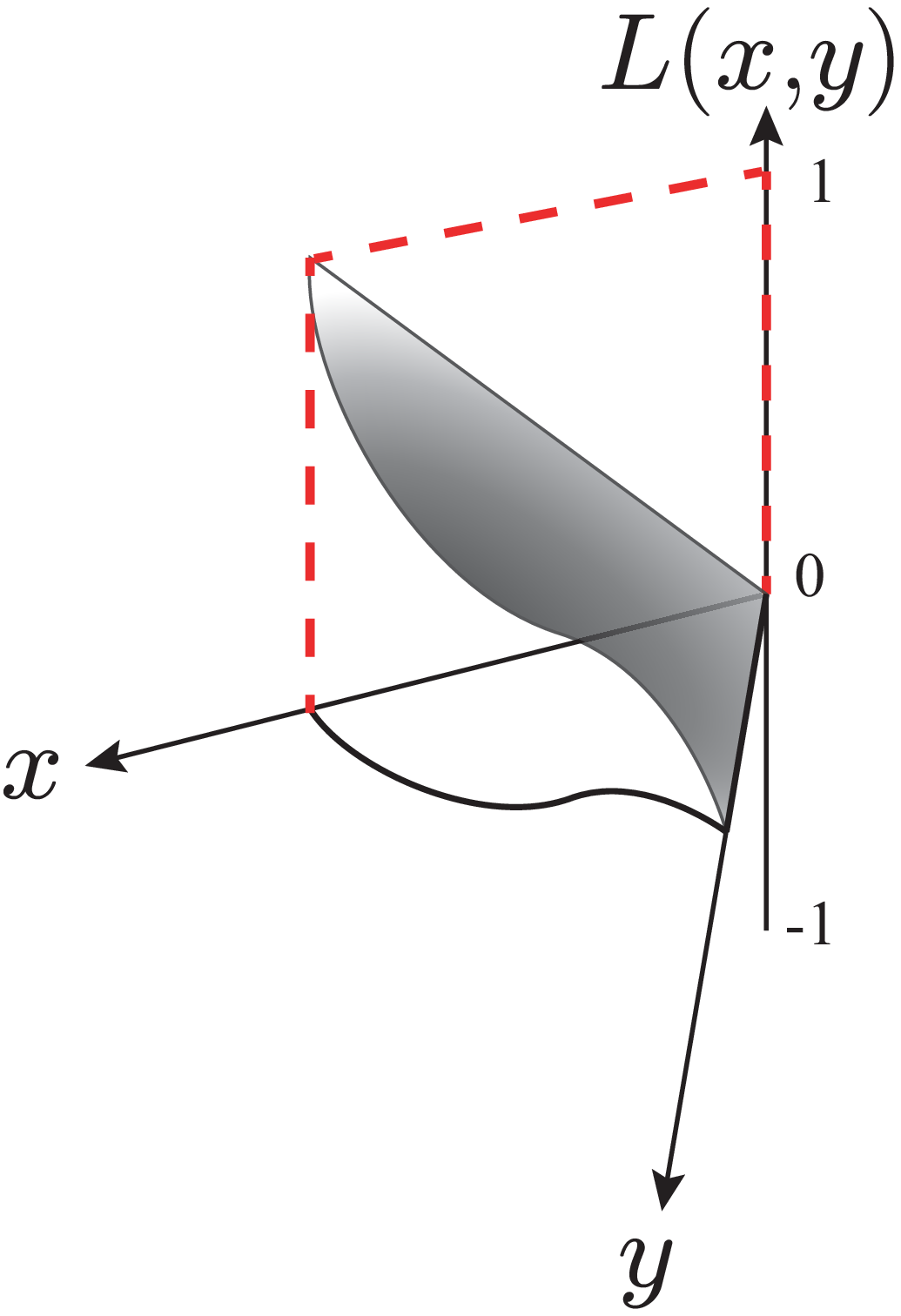}
    }
  \subfigure[Shape function for an interior node]{
    \includegraphics[width=6cm]{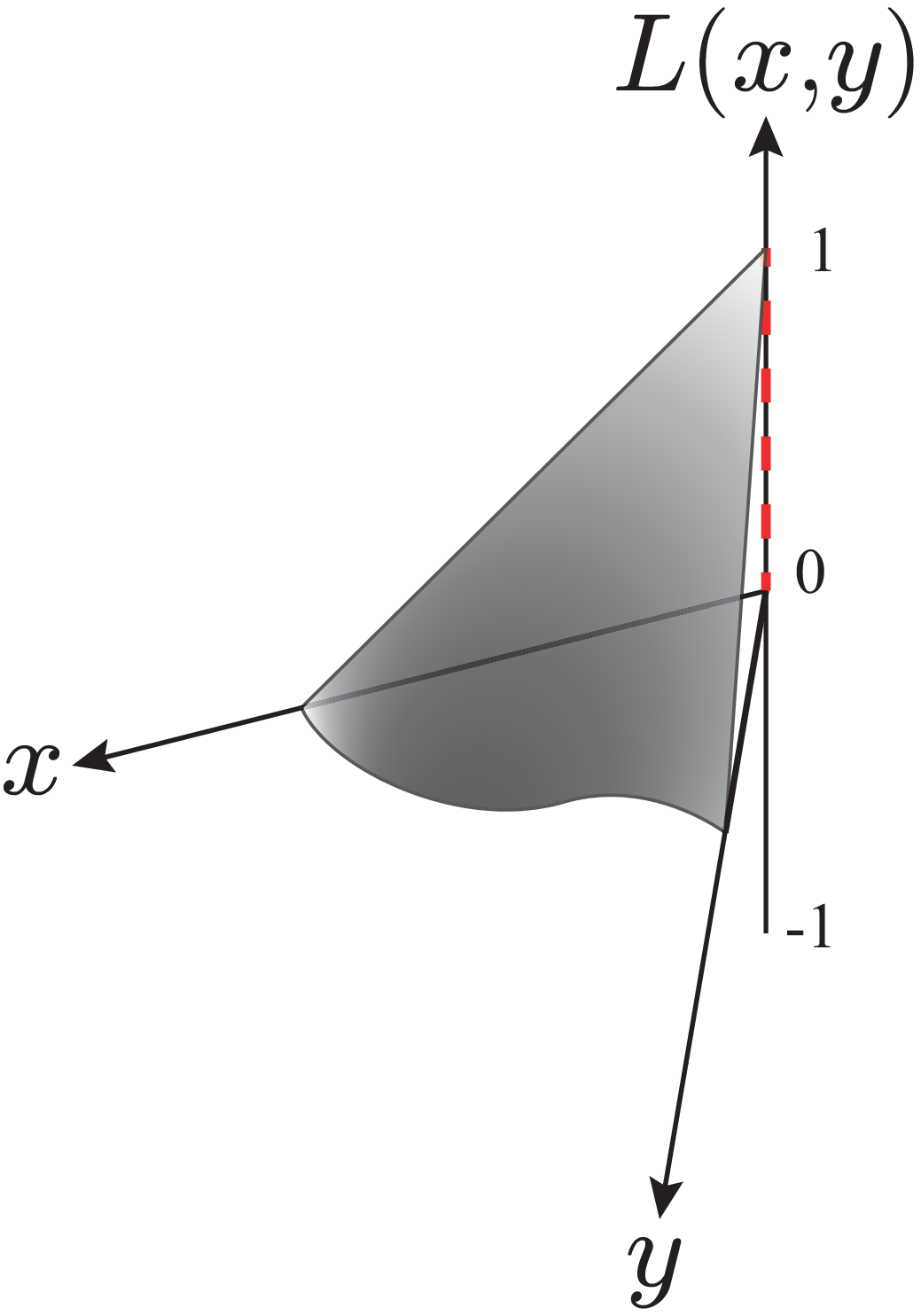}
  }
\caption{Shape functions for elements with curved edges. In (a) a shape function for a boundary node is exemplified. In (b), we see a shape function for the interior node. Note that it is zero along the curved boundary edge.}
\label{pfemshapefunctions}
\end{figure}

The reason for computing the shape functions on a reference element --- and thereby accepting the complication of distorted polynomial functions (cf. Figure~\ref{pfemshapefunctions}) --- lies in the fact that for Cartesian FEM, also the shape functions of the interior node are non-zero along the boundary edge/face. As a consequence, interior nodes contribute to boundary integrals. In the context of Dirichlet boundaries and boundary integrals, this shape function definition, however, leads to a grave disadvantage: the values along the boundary edge are always affected by nodal values of interior nodes. Since both Dirichlet boundary conditions and the boundary integral connected to the load transfer in FSI are of utter importance to the scenarios considered in this work, it was imperative to modify the shape function definition such that shape functions of inner nodes are zero along the curved edge. The shape function definition --- inspired by the so-called "p-FEM" in \cite{Sevilla2008a} --- is performed on the reference element. It is transformed to the global elements using a non-linear mapping ${\boldsymbol{\Phi}}$, which includes the NURBS definition. The mapping --- termed Triangle-Rectangle-Triangle (TRT) mapping --- is given as:

\begin{align}
    {\boldsymbol{\Phi}}(s,r) = (1-s-r) {\bf x}_2 + (s+r) {\bf C} \left( \frac{ {\Theta}_1 \; s   + {\Theta}_3 \; r }{s + r} \right)
\label{2DTRT}
\end{align}
  
Here, $s$ and $r$ denote the coordinates in the reference triangle. ${\bf x}_2$ is the global coordinate of the iterior node and ${\Theta}_{1,3}$ are the NURBS coordinates of the FE boundary nodes.

\begin{figure}[h]
\centering
 \includegraphics[scale=0.2,keepaspectratio=true]{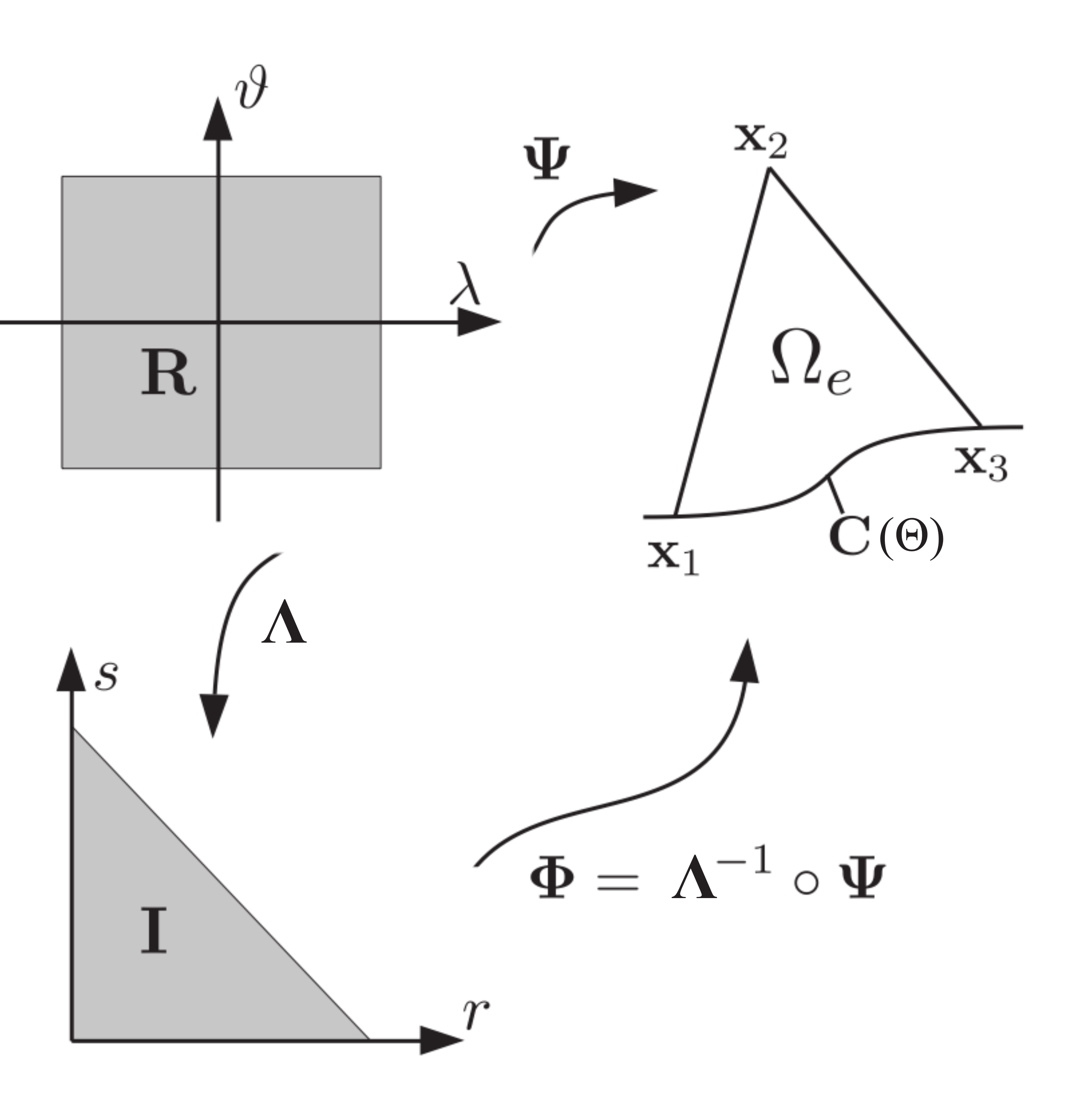}
 \caption{The TRT mapping from the reference triangle to the global triangle with one curved edge was derived from the original mapping from the reference bi-unit square to the curved triangle. By incorporating the bi-unit square, it can still be ensured that the NURBS direction and the interior direction are clearly separated, thus leading to straight interior edges even though the boundary edge is curved. ${\boldsymbol{\Psi}}$ is the mapping utilized in \cite{Sevilla2008a}, ${\boldsymbol{\Phi}}$ is the new TRT mapping.}
 \label{TRTmapping}
\end{figure}

\begin{figure}[h!]
\centering
    \includegraphics[width=7.2cm]{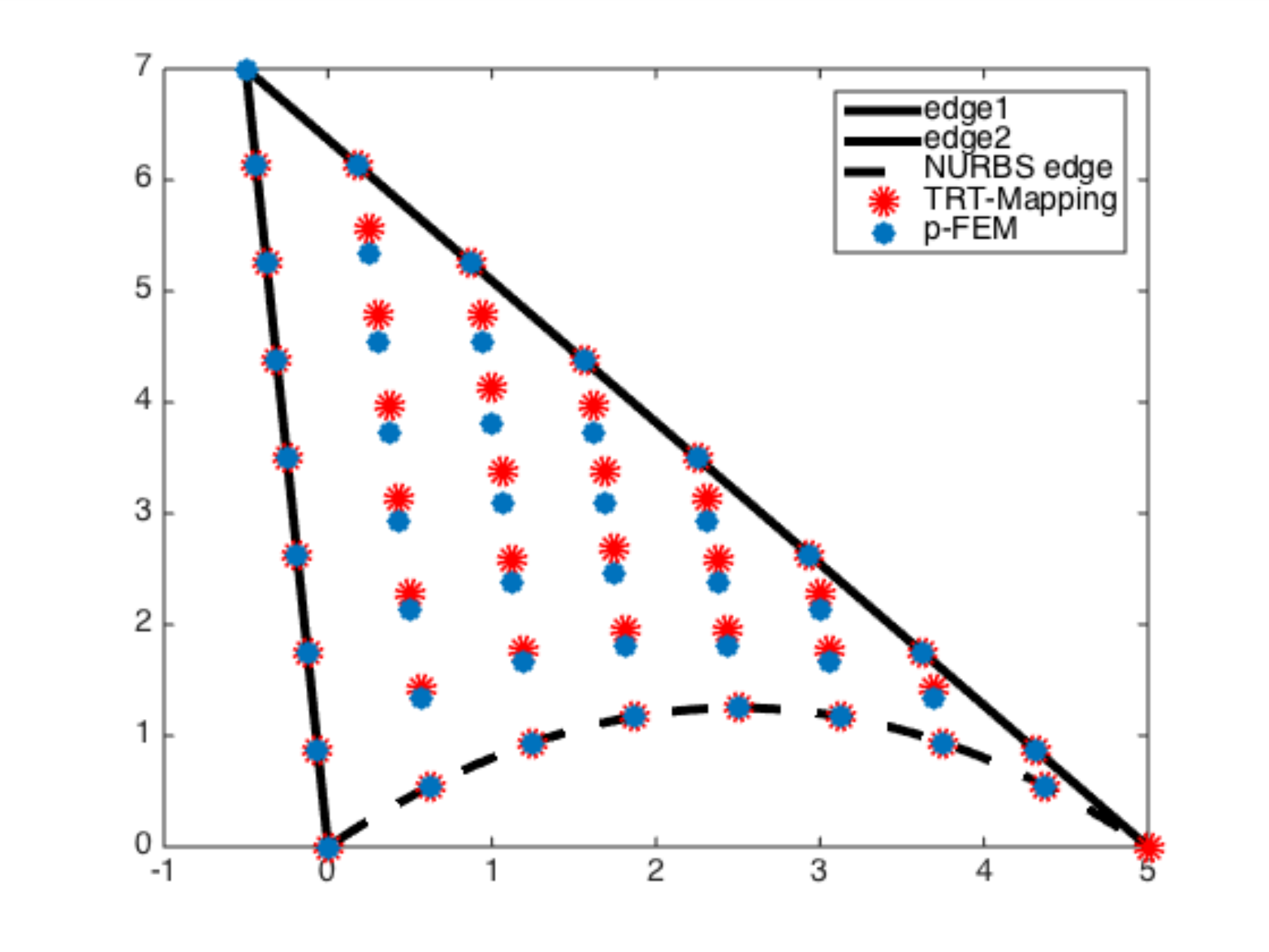}
\caption{Comparison of the quadrature point placement for the p-FEM and the TRT mapping for an element with one circular edge}
\label{pfemnefemcomparison1}
\end{figure}

 The derivation of the mapping is illustrated in Figure \ref{TRTmapping}. It is used for both the definition of the shape function and the placement of the quadrature points.
 Figure~\ref{pfemnefemcomparison1} compares the p-FEM mapping from \cite{Sevilla2008a}  with the mapping in Equation~\eqref{2DTRT}. In principle, the two mappings lead to a similar --- although still different --- distribution of quadrature points. Note however, that the p-FEM mapping has a singularity at one of the boundary nodes; thus excluding this point in case of a boundary integral. The TRT mapping has its singularity at the interior node.
 
\section{Numerical Methods: Coupling Approach} \label{sec-coupling}

In this work, we adopt a partitioned coupling approach, where the fluid and the structure are treated as individual fields and solved separately. The coupling conditions are explicitly incorporated as a means for interchanging information between the structure and the fluid.  Due to the use of individual solvers, one will usually encounter discrepancies between the discretizations of the structure and the fluid: starting from different levels of refinement, through different interpolation orders, to possibly even completely different discretization techniques. As a consequence, there is a demand for a strategy to interchange information/quantities between the individual fields. The most common approach is the "Neumann/Dirichlet" load transfer \cite{Braun2007}: Forces resulting from the fluid boundary stresses are projected onto the structure as a Neumann boundary condition, while the structural deformations are transferred to the fluid as a Dirichlet boundary condition.
 
Two main aspects have to be considered:

\begin{itemize}
\item temporal coupling: synchronization of individual fields,
\item spatial coupling: transfer of loads and deformation.
\end{itemize}

\subsection{Temporal Coupling}

For the temporal coupling, we employ a strong coupling approach. Its advantage is that
\--- in contrast to the often utilized weak coupling \--- the coupling conditions are
fulfilled after each time step. This is achieved via fixed-point iterations between the
structure and fluid within one time step until convergence. 

\begin{figure}[h]
\centering
 \includegraphics[scale=0.1,keepaspectratio=true]{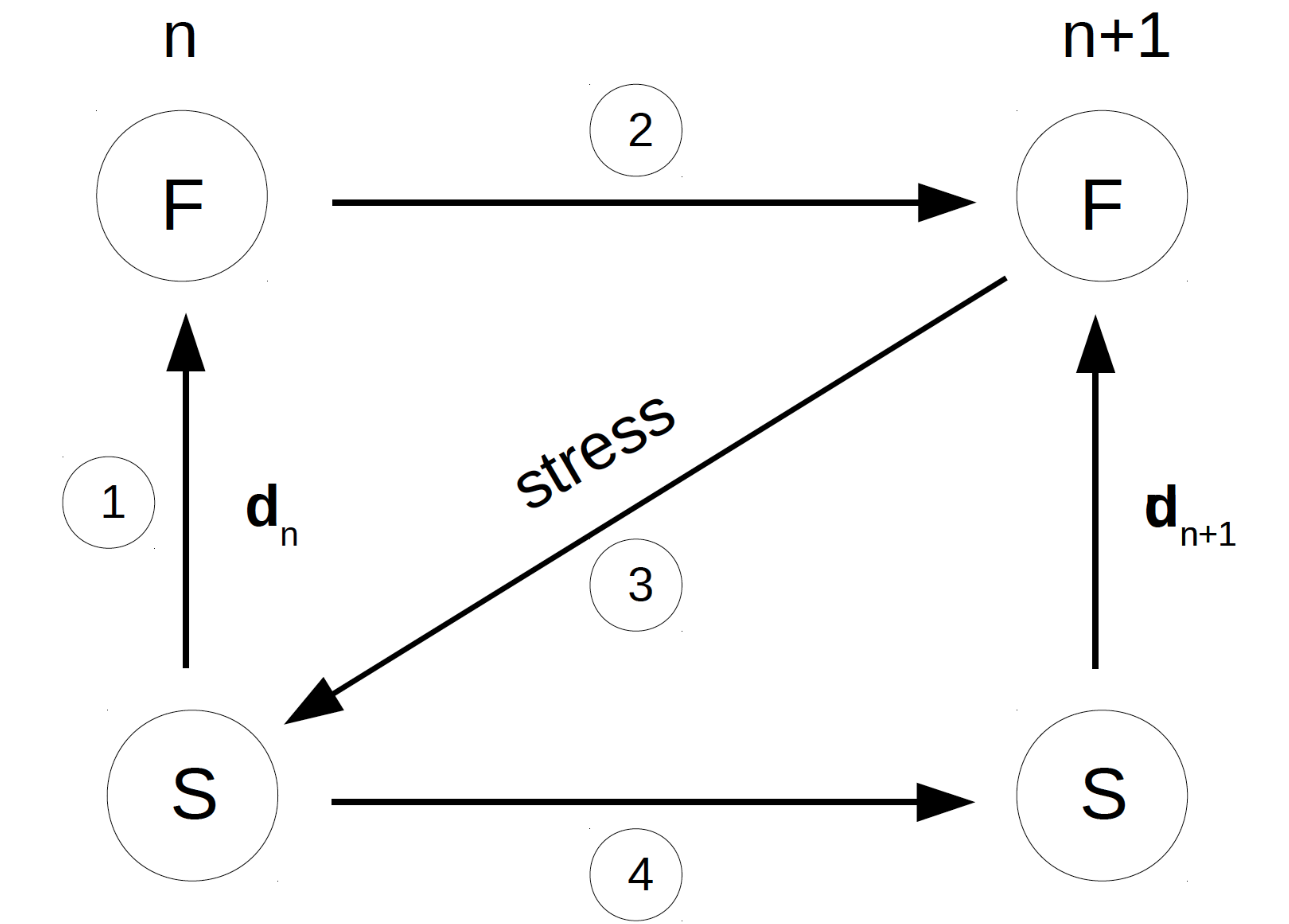}
 \caption{Schematic description of a strong coupling.}
 \label{strongcoupling}
\end{figure}

A schematic description of the coupling strategy is given in the following (see Fig. \ref{strongcoupling}):

\begin{enumerate}
\item compute an initial guess of the deformation $\tilde{\bf{d}}^{n+1}_k$ based on given quantities at $t_n$ and transfer it to the fluid,
\item deform the fluid domain based on $\tilde{\bf{d}}^{n+1}_k$ and calculate new fluid solution ${\bf{u}}^{n+1}_k$, $p^{n+1}_k$,
\item transfer stresses based on  ${\bf{u}}^{n+1}_k$, $p^{n+1}_k$ onto structure,
\item use the stresses to compute ${\bf{d}}^{n+1}_{k+1}$,
\item test convergence $ \lVert {\bf{d}}^{n+1}_{k+1} - {\bf{d}}^{n+1}_{k} \rVert \textless \epsilon $. If not converged go back to 2 and compute ${\bf{u}}^{n+1}_{k+1}$, $p^{n+1}_{k+1}$. 
\end{enumerate}

\subsection{Spatial Coupling}\label{sec:spatialcoupling}

For spatial coupling, we utilize a partitioned approach with individual solvers for both the structure and the fluid.  As a consequence, information has to be transferred from one solver to the other, which then enters the respective other simulation as a boundary condition: From structure to fluid these are the stresses, from fluid to structure these are deformations. In this work, we have chosen two possible approaches for spatial coupling: (1) the finite interpolation method, which is also applicable to standard FEM, and (2) direct integration of the discrete forces.

\subsubsection{Finite Interpolation Method} \label{sec:fim}

The general idea of this method is to transfer the discrete force at every fluid node onto its corresponding base point on the structure, and subsequently use the local non-zero basis function of the structural discretization to distribute the load to the corresponding nodes of the structure. In order to transfer the deformation back to the fluid domain, this procedure can be reversed. The individual steps are detailed in the following: \\

\noindent {\em 1) Calculation of discrete forces:} \\
The method requires discrete forces at each fluid node at the interface $\Gamma_{FS}$. Therefore, the discrete stresses calculated by the fluid solution have to be integrated over the interface faces and distributed to each node. Given the discrete stress $\boldsymbol{\sigma}^f_i$ at each node on the boundary, the stress distribution is defined as follows:

\begin{equation}
 \boldsymbol{\sigma} \left( r \right) = \sum_i L_i \left( r \right) \boldsymbol{\sigma}^f_i,
\end{equation}

\noindent where $L_i$ denotes the Lagrangian shape functions and $r$ is the parametric coordinate along the interface. In order to regain the discrete force at each individual node, the stress has to be multiplied by the corresponding shape function and then integrated, resulting in the formula

 \begin{equation}\label{discreteforcenode}
 {\bf{F}}_i^f = \int_{\Gamma_{FS}} L_i  \left( r \right) \sum_j \left[ L_j \left( r \right) \boldsymbol{\sigma}^f_j \cdot {\bf{n}} \right]  \; dr,
\end{equation}

\noindent where $\bf{n}$ is the normal vector at the interface. \\

\noindent {\em 2) Structural base point computation: }\\
In the next step, the base point on the structure of each fluid interface node has to be found. With our interface representation, this is straightforward: Since the spline representations are identical for both fluid and structure, the base point is identified through the local spline coordinate $\Theta_j$ associated with the specific NEFEM node. During preprocessing, it can, e.g., be calculated by using an orthogonal projection of a point onto a spline based on \cite{Hu2005}. For multiple patches, the corresponding patch has to be identified before the parametric coordinate is computed. \\

\noindent {\em 3) Transfer of forces onto the structure:} \\
We want to obtain the discrete force at every single control point of the spline. Using its NURBS basis function $R(\Theta)$, we can simply sum over all discrete forces multiplied by the basis function evaluated at the corresponding local parametric coordinate with: 

\begin{equation}
{\bf{F}}_i^s = \sum_j^{nn^f_{\Gamma_{FS}}} R_i \left( \Theta _j \right) {\bf{F}}_j^f.
\end{equation}

\noindent Here, ${nn^f_{\Gamma_{FS}}}$ stands for the number of fluid nodes on the interface. The NURBS basis functions fulfill the partition of unity property, meaning that

\begin{equation}
\sum_i^{n_{CP}} R_i \left( \Theta \right) = 1 \quad \forall \; \Theta,
\end{equation}

\noindent holds. Therefore, the overall force acting on the interface is conserved during the projection. \\

\noindent {\em 4) Transfer of deformation onto the fluid: }\\
The computed structural deformation has to be transferred back to the fluid in order to account for the mesh deformation. Once again, this procedure is straightforward due to the use of IGA. The deformation is computed at every spline control point. Knowing the parametric coordinates, we have to sum up the contributions of all control points that have non-zero NURBS basis function at that point. This is done in the following way:

\begin{equation}
{\bf{d}}_i^f = \sum_j^{n_{CP}} R_j \left( \Theta_i \right) {\bf{d}}_j^s.
\end{equation}

Here, ${\bf{d}}_j^s$ is the computed displacement on the structural side --- a displacement applied to the control points ---, whereas ${\bf{d}}_i^f$ is the displacement of fluid node $i$ --- a displacement applied on FE node level.

\subsubsection{Direct Transfer} \label{sec:directtransfer}

If NEFEM is applied on the fluid side and IGA on the structural side, the geometry of the
coupling interface is identical on both sides. Thus, a direct integration can be applied to obtain the fluid forces in sense of a weighted residual method \cite{deBoer2008}.

\noindent {\em 1) Calculation of discrete forces:}

Using the NURBS basis functions as test functions, the right-hand-side of the structural problem can be formulated as follows:

\begin{equation}\label{NEFEMdiscreteforcenode1}
 {\bf{F}}_i^s = \int_{\Gamma_{FS}} R_i(\Theta) \left(\boldsymbol{\sigma}^f(\Theta) \cdot {\bf{n}}(\Theta) \right)  \; d\Theta.
\end{equation}

The fluid stresses $\boldsymbol{\sigma}^f$ can be computed from the fluid solution and the face normals from the fluid grid. For NEFEM, the face normal can be determined exactly by evaluating the NURBS. The boundary integral can be compute piecewise by computing it on every fluid grid face:

\begin{equation}\label{NEFEMdiscreteforcenode2}
 {\bf{F}}_i^s =  \sum^{n_F}_e \int_{\Gamma_{FS}^e} R_i(\Theta) \left(\boldsymbol{\sigma}^f(\Theta) \cdot {\bf{n}}(\Theta) \right)  \; d\Theta.
\end{equation}

By applying Gaussian integration with $n_{GP}$ being the number of Gauss points we obtain

\begin{equation}\label{NEFEMdiscreteforcenode3}
 {\bf{F}}_i^s =  \sum^{n_F}_e \sum^{n_{GP}}_g w(\Theta_g) R_i(\Theta_g) \left(\boldsymbol{\sigma}^f(\Theta_g) \cdot {\bf{n}}(\Theta_g) \right).
\end{equation}

If we insert now the Lagrangian representation of the fluid solution, the forces acting on the right-hand-side of the structural problem can be evaluated as
\begin{equation}\label{NEFEMdiscreteforcenode4}
 {\bf{F}}_i^s =  \sum^{n_F}_e \sum^{n_{GP}}_g w(\Theta_g) R_i(\Theta_g) \left(\sum_j L_j \boldsymbol{\sigma}^f_j(\Theta_g) \cdot {\bf{n}}(\Theta_g) \right).
\end{equation}

Following this idea, the forces are already available on every control point and can be used as a right-hand-side for the system of equations of the structural problem. An additional projection method is not needed anymore.

Steps {\em 2)} and {\em 4)} remain the same as in section \ref{sec:fim}. Because the numerical integration is conducted on the same geometry, the resulting method is consistent and conservative \cite{deBoer2008}.

\section{Numerical Examples} \label{sec-testcases}

In order to evaluate the performance of the new fluid-structure interaction approach, we
have chosen three different test cases. The first was proposed by Sch\"afer et al.
~\cite{featflow}, the second one by Wall~\cite{Wall99} and the third one by Turek et al.~\cite{turek_proposal_2007}. The first two test cases involve steady laminar incompressible channel flow around a rigid and elastic cylinder respectively. The third test case is unsteady lamir flow around a structure that oscillates due to flow-induced vibrations. Thereby, the adapted NEFEM method with direct coupling can be tested for simple problems. What we expect to see is that the convergence rate will remain the same as for standard FEM, as we still utilize the same shape function definition, but the absolute value of the error will decrease. The latter is due to a decrease in geometrical error due to integration over a smooth and ---- in the rigid case --- even exact domain.

\subsection{Cylinder}
\label{subsec-cylinder}

\begin{figure}[h]
\centering
 \resizebox{.7\linewidth}{!}{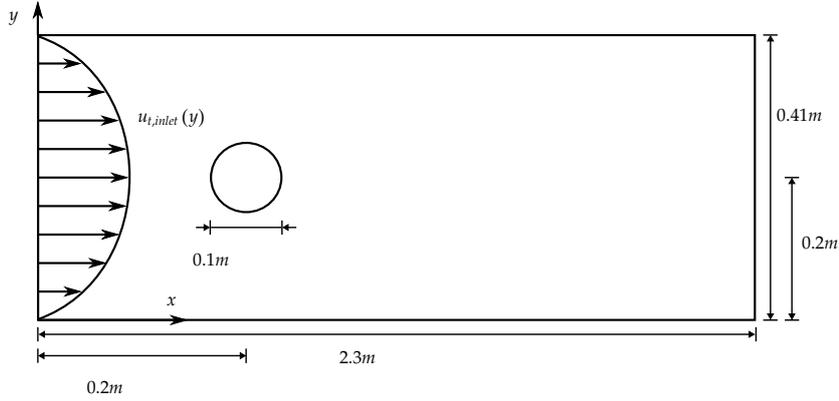}
 \caption{Sketch of the rigid featflow benchmark.}
 \label{turekbenchmark}
\end{figure}

The first test case is used to compare the influence of using NEFEM in contrast to FEM and involves laminar flow around an attached, rigid cylinder. The configuration is proposed on the featflow benchmark site \cite{featflow}. In particular this is the ``DFG flow a around cylinder benchmark 2D-1''. A sketch of the computational domain is provided in Figure \ref{turekbenchmark}.

A parabolic inflow velocity is given at the inlet:

\begin{equation} \label{eq:para_inflow}
{\bf u} \left( 0, y, t \right)  = \left( \frac{4 U y \left( H - y \right)}{H^2} , 0 \right).
\end{equation}

No-slip conditions are applied to the upper and lower wall as well as to the cylinder and a free-flow boundary condition to the outlet. 
All parameters used for the simulation can be found in Table \ref{table:rigidfeatflow}. The Reynolds number is defined as

\begin{equation}
Re = \frac{\rho \overline{U} D}{\mu},
\end{equation}

where $\overline{U}$ is the average inflow velocity.

\begin{table}[h]
\begin{minipage}[t]{0.5\textwidth}
\centering
\begin{tabular}{|l|c|c|}
\hline
 parameter & identifier & value\\
 \hline
inflow velocity & U & 0.2*1.5 $m/s$\\
dynamic viscosity & $\mu$ & 0.001 $kg/m s$ \\
density & $\rho$ & 1.0 $kg/m^3$ \\
cylinder diameter & D & 0.1 $m$ \\
Reynolds number & Re & 20\\
\hline
\end{tabular} \caption{Rigid featflow benchmark parameters.}
\label{table:rigidfeatflow}
\end{minipage}
\begin{minipage}[t]{0.5\textwidth}
\centering
\begin{tabular}{|l|c|c|}
\hline
  & elem. & elem. on cylinder\\
 \hline
Grid 1 &  5411 & 50 \\
Grid 2 & 13863 &  100  \\
Grid 3 & 55442 &  200  \\
Grid 4 & 221808 & 400  \\
Grid 5 & 887232 & 800 \\
\hline
\end{tabular} \caption{Grids used for the featflow benchmark simulations.}
\label{parametergrid}
\end{minipage}
\end{table}

For the given Reynolds number, the flow is steady. Therefore, the steady Navier-Stokes equations are solved. 

An unstructured grid with 50 boundary elements on the cylinder and 5411 elements in total is used as a starting point. Four further grids are considered, which are obtained by a simple h-refinement of the initial grid. This is done by splitting every edge in half. An overview over the utilized grids is given in Table \ref{parametergrid}. In case of NEFEM, the geometry of the cylinder is represented by a second order NURBS line with 10 control points.

\begin{figure}[h!]
    \centering
    \large
    \resizebox{.45\linewidth}{!}{\input{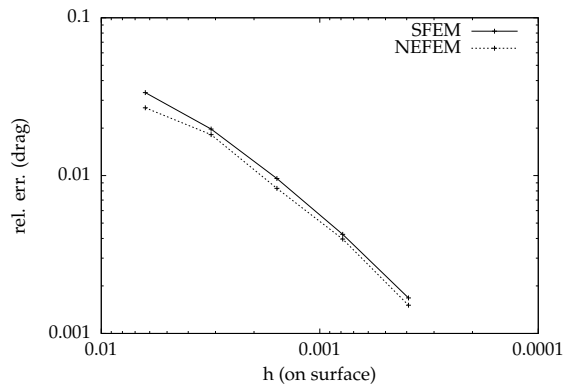}}
    \caption{Grid convergence for non-moving clamped cylinder}
    \label{convergenceSpringrigid}
\end{figure}

The reference value that serves as a basis for the comparison between the results for NEFEM and standard FEM is the drag coefficient. In both cases --- NEFEM and FEM --- the drag coefficient computed on the final grid (NEFEM: $c_d = 5.570928 $ and FEM: $c_d = 5.569983$) is in good agreement with the reference value of $c_d = 5.579535$ given by featflow. We therefore consider both implementations to be valid for this test case. In addition, Figure \ref{convergenceSpringrigid} gives a convergence plot for the drag obtained with the two methods. The featflow result is employed as a reference result. As expected, both methods feature the same convergence rate. However, one can observe that NEFEM slightly improves the error constant.

\subsection{Flexible Cylinder} \label{subsec-breathing}

\begin{figure}[h]
    \centering
    \resizebox{.7\linewidth}{!}{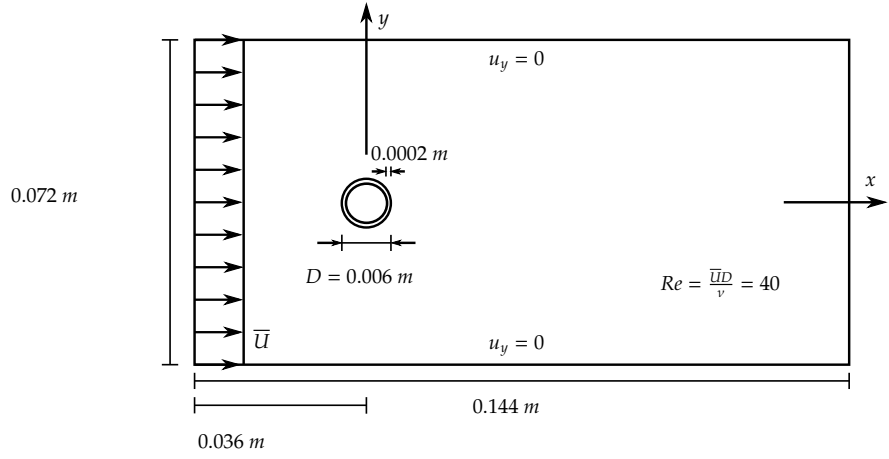}
    \caption{Sketch of the flexible cylinder test case.}
    \label{domainBreathCyl}
\end{figure}

The second test case can be viewed as an extention of test case \ref{subsec-cylinder}. Instead of a rigid cylinder, it now considers an elastic structure \cite{Wall99}.  The dimensions of the domain are detailed in Fig. \ref{domainBreathCyl}. The inflow velocity is constant; on the lower and upper wall a slip condition is applied. A no-slip boundary condition is used on the cylinder. The mass enclosed by the cylinder is neglected. The elastic structure is completely fixed on its rightmost point (on the horizontal symmetry axis). In order to avoid asymmetric results induced by numerical errors, the displacement of the leftmost cylinder point is suppressed in y-direction.

\begin{table}[h!]
\captionsetup{justification=raggedright} 
\centering
\begin{minipage}[b]{0.5\linewidth}
\centering
\begin{tabular}{|l|c|c|}
 \hline
 parameter & identifier & value\\
 \hline
inflow velocity & U & 0.10282 $m/s$\\
dynamic viscosity & $\mu$ & 0.0000182 $kg/m s$ \\
density (fluid) & $\rho^f$ & 1.18 $kg/m^3$ \\
cylinder diameter & D & 0.006 $m$ \\
Reynolds number & Re & 40\\
 \hline
Young's modulus& $E$& 100 Pa\\
Poisson number&$\nu^s$&0.3\\
density (structure)&$\rho^s$& 1000 $kg/m^3$\\
\hline
\end{tabular} 
\end{minipage}%
\begin{minipage}[t]{0.5\textwidth}
\centering
\begin{tabular}{|l|c|c|}
\hline
  & elem. & elem. on cylinder\\
 \hline
Grid 1 &    2352 &   40 \\ 
Grid 2 &    9408 &   80 \\
Grid 3 &   37632 &  160 \\
Grid 4 &  150528 &  320 \\
Grid 5 &  602112 &  640 \\
Grid 6 & 9633792 & 2560 \\
\hline
\end{tabular} 
\end{minipage}
\par
\begin{minipage}[b]{0.5\textwidth}
\caption{Flexible cylinder parameters.}
\label{table:breathCyl}
\end{minipage}%
\begin{minipage}[b]{0.5\textwidth}
\caption{Grids used for the flexible cylinder simulations.}
\label{parametergridbreathCyl}
\end{minipage}

\end{table}

\begin{figure}[h!]
    \centering
	\subfigure[Undeformed cylinder]{
	\includegraphics[trim=0.1cm 0.1cm 6.0cm 1cm,clip,height=0.4\linewidth]{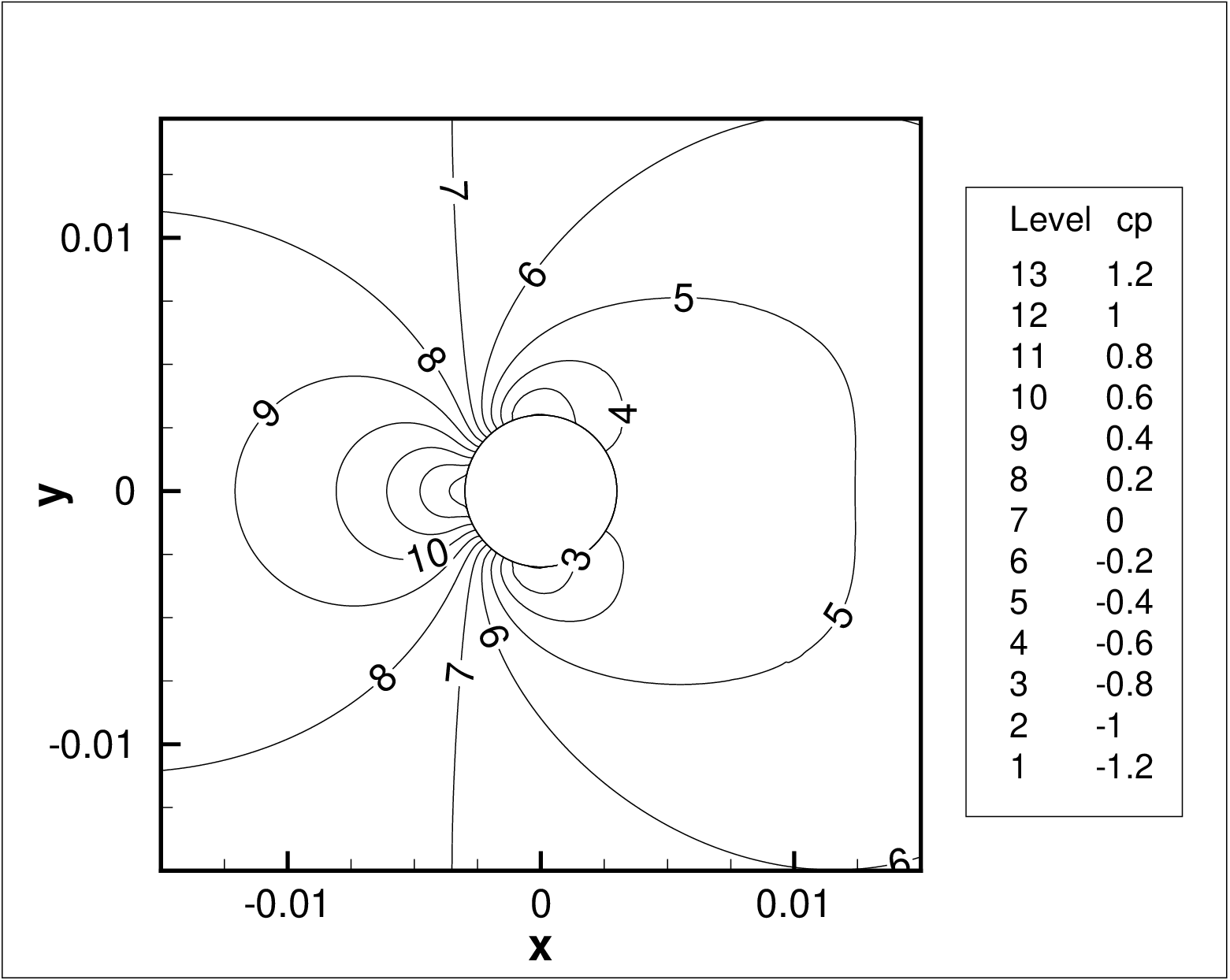}
	}
	\hfill
	\subfigure[Deformed cylinder]{
	\includegraphics[trim=0.1cm 0.1cm 0.5cm 1cm,clip,height=0.4\linewidth]{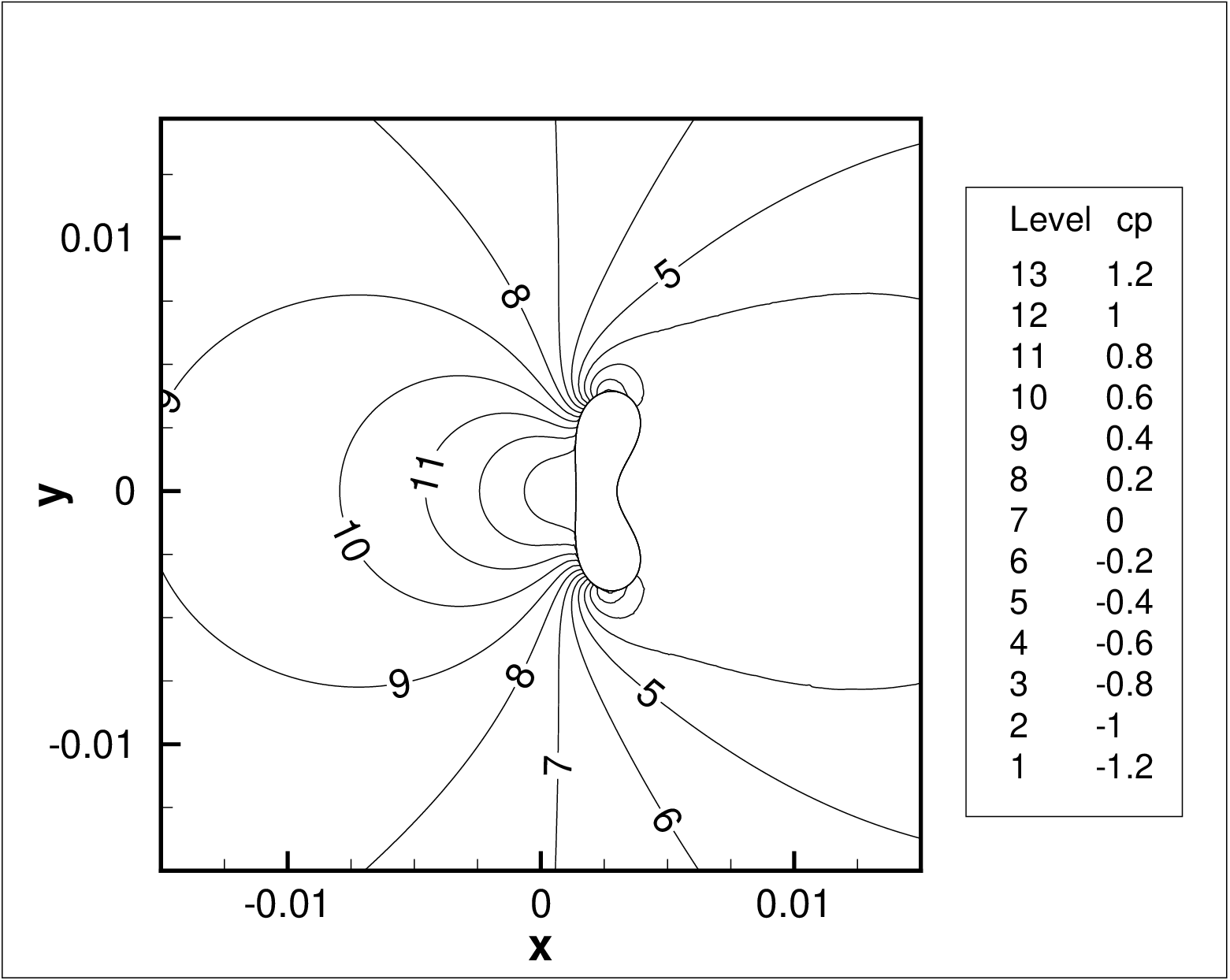}
	}
    \caption{Pressure distribution for the undeformed and the deformed cylinder.}
    \label{Pressure}
\end{figure}

The given Reynolds number of 40 leads to a steady solution and a significant deformation of the cylinder. Pressure distributions around the rigid and deformed cylinder are given in Fig. \ref{Pressure}.

On the fluid side,  an unstructured grid with 40 boundary elements on the cylinder and 2352 elements in total is generated initially. Similar to the first test case, six additional grids are created by applying simple h-refinement to the initial grid. An overview of the grids used is given in Table \ref{parametergridbreathCyl}. The geometry of the structure is defined by a second-order NURBS with 40 elements in circular direction and 5 elements in radial direction. Because the effect of NEFEM on FSI should be isolated, the discretization of the elastic structure remains unchanged within the presented study. 

\begin{figure}[h!]
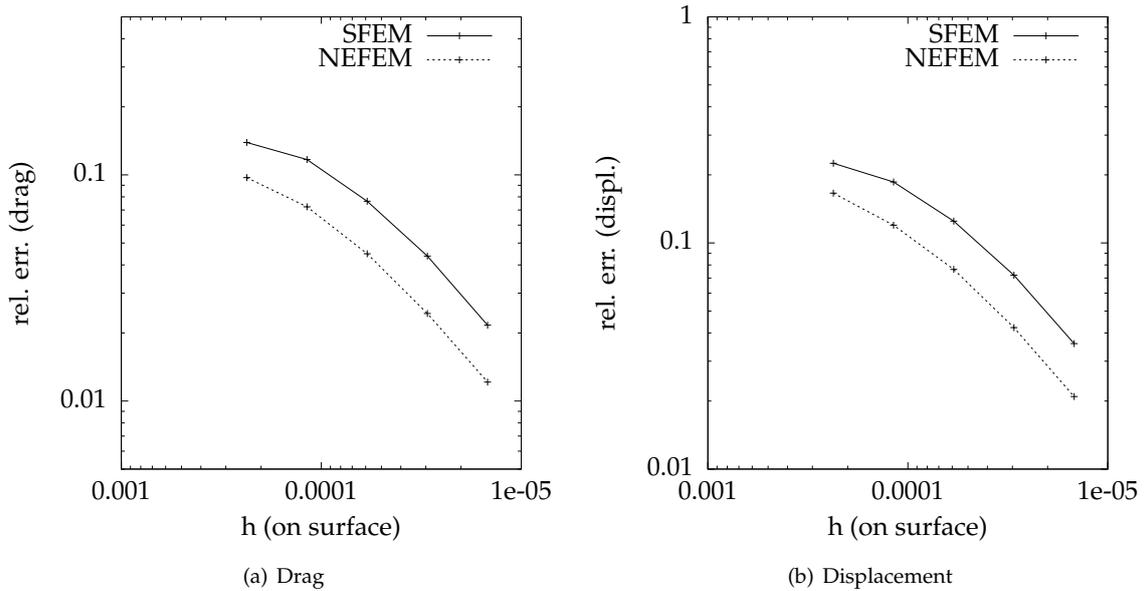

    \centering
    \large
    \subfigure[Drag]{
	    \resizebox{.45\linewidth}{!}{\input{./results/breathingCyl/ConvergencyDrag.tex}}
    }
    \subfigure[Displacement]{
    \resizebox{.45\linewidth}{!}{\input{./results/breathingCyl/ConvergencyDisplacement.tex}}
    }
    \caption{Grid convergence for the flexible cylinder.}
    \label{conBreathCylDrag}
\end{figure}

Two reference quantities are used to assess the convergence: (1) drag on the cylinder and (2) the displacement of the leftmost cylinder point. The results are presented as relative results to grid 6 in Fig. \ref{conBreathCylDrag}. As in the rigid case, the slope of the graphs are almost identical, but the relative error of NEFEM is clearly lower. In the presented case, the offset between both graphs is in the order of one refinement level.   

\subsection{Moving Flag}

\begin{figure}[h!]
\centering
    \includegraphics[width=.8\linewidth]{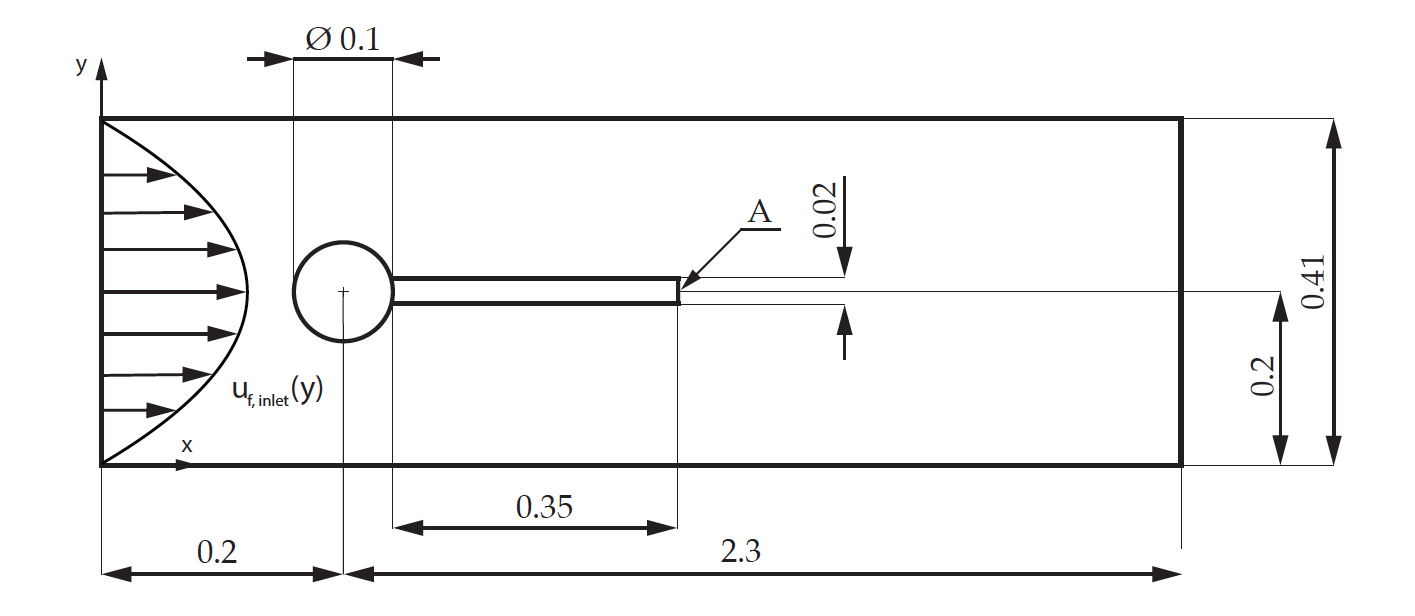}
    \caption{Sketch of the moving flag featflow benchmark.}
    \label{domainMovFlag}
\end{figure}

The third testcase is used in order to show the functioning of the method for transient
FSI problems. The configuration is proposed in~\cite{turek_proposal_2007}
and is referred to as "FSI benchmark - FSI2 ". A sketch of the computational domain is
given in Fig. \ref{domainMovFlag}. The deformable elastic structure, or "flag", is fixed at the flag-cylinder interface.

\begin{table}[h!]
\centering
\begin{minipage}[b]{0.5\textwidth}
\begin{tabular}{|l|c|c|}
\hline
 parameter & identifier & value\\
 \hline
inflow velocity & U & 1.0*1.5 $m/s$\\
dynamic viscosity & $\mu$ & 1.0 $kg/m s$ \\
density & $\rho^f$ & 1000.0 $kg/m^3$ \\
cylinder diameter & D & 0.1 $m$ \\
Reynolds number & Re & 100\\
\hline
Young's modulus& $E$& 1.4 10e6  $Pa$\\
Poisson number&$\nu^s$&0.4\\
density (structure)&$\rho^s$& 10000.0 $kg/m^3$\\
\hline
\end{tabular} 
\end{minipage}%
\begin{minipage}[b]{0.5\textwidth}
\centering
\begin{tabular}{|l|c|c|}
\hline
  & elements & elements on flag \\
 \hline
Grid 1 &    3572 &   75 \\ 
Grid 2 &    14288 &   150 \\
Grid 3 &   57152 &  300 \\
Reference & 162134 & 584 \\
\hline
\end{tabular} 
\end{minipage}
\par
\begin{minipage}[b]{0.5\textwidth}
\caption{Moving flag featflow benchmark parameters.}
\label{table:movingfeatflow}
\end{minipage}%
\begin{minipage}[b]{0.5\textwidth}
\caption{Grids used for moving flag simulations.}
\label{parametergridmovFlag}
\end{minipage}

\end{table}

\begin{figure}[h!]
    \centering
    \large
    \includegraphics[trim=0.1cm 0.1cm 0.5cm 0.1cm,clip,width=\linewidth]{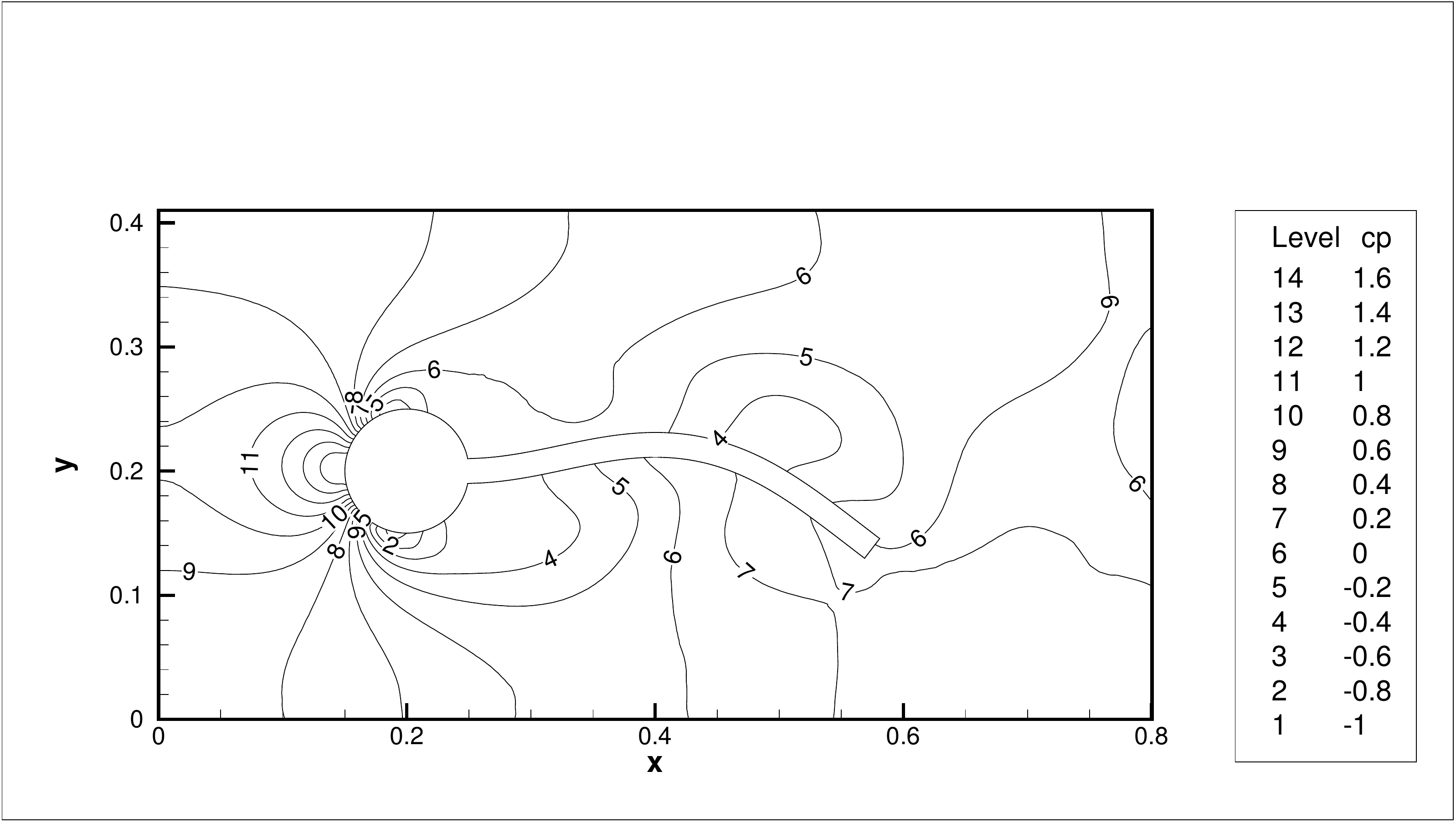}
    \caption{Pressure distribution for deformed flag at time t=13.6s}
   \label{fig:presFlag}
\end{figure}

Similar to Section \ref{subsec-cylinder} a parabolic velocity profile is set at the inlet.
No-slip conditions are applied to the upper and lower wall as well as to the cylinder and
the flag and a free-flow boundary conditions is applied to the outlet. All parameters
relevant for the simulations are shown in Table \ref{table:movingfeatflow}. For the given
Reynolds number of 100 the interaction of the flow and the elastic structure causes
self-induced oscillations of the structure which results in a quasi-steady periodic wave-like deformation of the flag. The pressure distribution at time t=13.6 s is given
in Fig.~\ref{fig:presFlag}. 

We discretize the fluid domain using an unstructured grid with 75 boundary elements on the
flag surface and 3572 elements in total. Similar to the previous test cases we create two
other grids by applying simple h-refinement to the intitial grid. As reference for the
convergence study we use the solution computed on a fine grid with 584 elements on the flag surface. An overview over the grids used is given in
Table \ref{table:movingfeatflow}. The structure is discretized using a
quadratic NURBS with 60x15 elements. It is kept constant throughout the study. In all
computations a time step of $\Delta t=0.002s$ was applied.

\begin{figure}[h!]
    \centering
    \large
    \resizebox{.45\linewidth}{!}{\input{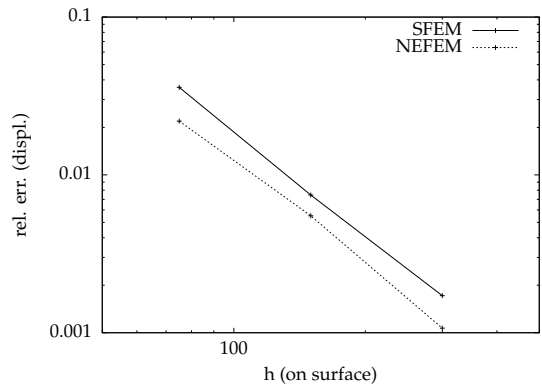}}
    \caption{Grid convergence of amplitude of vertical displacement at point A.}
    \label{fig:convergenceFlag}
\end{figure}

The quantity for comparison is the amplitude of the vertical displacement of point A,
located at the tail of the flag. Fig.~\ref{fig:convergenceFlag} shows the results as
relative error to the reference solution. As in the previous test cases, the offset between both graphs is identical, but again the relative error of NEFEM is clearly lower. 

\newpage

\section{Conclusion} \label{sec-conclusion}

The paper presents a novel coupling scheme for fluid-structure-interaction problems.
Building on the recent advances in structural mechanics based on isogeometric analysis,
the idea is to involve the splines also in the flow simulations. This was achieved by
employing the NURBS-Enhanced Finite Element Method for the fluid. A new geometrical
mapping from physical to reference space ensuring that Dirichlet boundary conditions are fulfilled on the spline-based
surface was presented for NEFEM. With that, a direct transfer of the necessary coupling variables is
possible, as the interface description is identical for both the structure and the fluid.
In addition to the simplified implementation of the coupling, an increased accuracy of the
flow solution --- and with this the FSI solution --- was to be expected due to the
matching interface representation. This expected improvement was confirmed with steady computations
on a rigid and deformable cylinder as well as transient
computations of flow-induced vibrations, see Section~\ref{sec-testcases}. Even if the convergence rate remains the same,
the resulting error of NEFEM compared to standard finite elements can be decreased. With
the second and third example it was shown that the method is directly applicable to elastic configurations using isogeometric analysis on the structural side. Finally, the presented method can be used to generate fluid grids for more complex problems with existing tools and allows for improved accuracy with minor implementation effort.

\section{Acknowledgements}
The authors gratefully acknowledge the support of the German research foundation DFG grant ``Geometrically Exact Methods for Fluid-Structure Interaction''. Furthermore, we thank Daniel Hilger for his support and ideas.

\bibliographystyle{elsarticle-num}
\section*{References}
\bibliography{references}

\end{document}